\documentclass[a4paper,12pt]{article}

\usepackage[margin=2.1cm]{geometry}
\usepackage{graphicx}
\usepackage{subfigure}
\usepackage[square,numbers]{natbib}
\usepackage{lipsum}
\usepackage{caption}
\usepackage{pdfpages}
\usepackage{amsmath}
\usepackage{multicol}
\usepackage{multirow}
\usepackage{float}
\usepackage{url}
\usepackage{appendix}
\usepackage{bm}
\usepackage{amssymb}

\usepackage{mwe}

\usepackage{soul, url}

\newcommand{\bfE}{\mathbf{E}}
\newcommand{\bfF}{\mathbf{F}}
\newcommand{\bfL}{\mathbf{L}}
\newcommand{\bfW}{\mathbf{W}}
\newcommand{\bfA}{\mathbf{A}}
\newcommand{\bfnu}{\boldsymbol{\nu}}
\newcommand{\bfchi}{\boldsymbol{\chi}}
\newcommand{\bfg}{\mathbf{g}}

\newcommand{\bfp}{\mathbf{p}}
\newcommand{\bfx}{\mathbf{x}}
\newcommand{\bfa}{\mathbf{a}}
\newcommand{\bfy}{\mathbf{y}}
\newcommand{\bfv}{\mathbf{v}}
\newcommand{\bfd}{\mathbf{d}}
\newcommand{\bfb}{\mathbf{b}}
\newcommand{\bfu}{\mathbf{u}}
\DeclareMathOperator{\curl}{curl}

\newcommand{\epsr}{\epsilon_{\mathrm{r}}}
\newcommand{\mur}{\mu_{\mathrm{r}}}

\newcommand{\hs}{h_{\mathrm{s}}}
\newcommand{\bfxi}{\boldsymbol{\xi}}
\newcommand{\hv}{h_{\mathrm{v}}}
\newcommand{\Nt}{N_{\mathrm{elements}}^{\mathrm{tetra}}}
\newcommand{\Nw}{N_{\mathrm{elements}}^{\mathrm{wedge}}}
\newcommand{\Nh}{N_{\mathrm{elements}}^{\mathrm{hexa}}}
\newcommand{\Nv}{N_{\mathrm{vertices}}}
\newcommand{\Ndof}{N_{\mathrm{DoF}}}

\begin{document}

\title{A High-Order Ultra-Weak Variational Formulation for Electromagnetic Waves Utilizing Curved Elements}

\author{Timo L\"ahivaara${}^{a}$, William F. Hall${}^{b}$, Matti Malinen${}^{c}$,\\ Dale Ota${}^{b}$, Vijaya Shankar${}^{b}$, and Peter Monk${}^{d}$\\
\normalsize ${}^{a}$Department of Technical Physics, University of Eastern Finland, Kuopio, Finland\\
\normalsize ${}^{b}$HyperComp Inc., Westlake Village, USA\\
\normalsize ${}^{c}$Kuava Ltd., Kuopio, Finland\\
\normalsize ${}^{d}$Department of Mathematical Sciences, University of Delaware, Newark, USA
}

\maketitle

\section*{Abstract}

The Ultra Weak Variational Formulation (UWVF) is a special Trefftz
discontinuous Galerkin method, here applied to the time-harmonic
Maxwell’s equations.  The method uses superpositions of plane waves to
represent solutions element-wise on a finite element mesh. We focus on
our parallel UWVF implementation, called \emph{ParMax}, emphasizing
high-order solutions in the presence of scatterers with piecewise
smooth boundaries. We explain the incorporation of curved surface
triangles into the UWVF, necessitating quadrature for system matrix
assembly. We also show how to implement a total field and scattered
field approach, together with the transmission conditions across an
interface to handle resistive sheets.  We note also that a wide
variety of element shapes can be used, that the elements can be large
compared to the wavelength of the radiation, and that a low memory
version is easy to implement (although computationally costly). Our
contributions are illustrated through numerical examples demonstrating
the efficiency enhancement achieved by curved elements in the
UWVF. The method accurately handles resistive screens, as well as
perfect electric conductor and penetrable scatterers. By employing
large curved elements and the low memory approach, we successfully
simulated X-band frequency scattering from an aircraft. These
innovations demonstrate the practicality of the UWVF for industrial
applications.

\newpage

\section{Introduction}

A Trefftz method~\cite{Trefftz26} for approximating a linear partial
differential equation is a numerical method using local solutions of
the underlying partial differential equation as basis functions.  The
version we shall study here, the Ultra Weak Variational Formulation
(UWVF) of Maxwell's equations, is a Trefftz type method for
approximating the solution of Maxwell's equations on a bounded domain
due to O. Cessenat and B. Despr\'es~\cite{cessenat_phd,cessenat03}.
The UWVF uses a finite element computational grid, classically
composed of tetrahedral elements, and plane wave solutions of
Maxwell's equations on each element.

A study of this method from the point of view of symmetric hyperbolic
systems was presented in \cite{Huttunen2007} where the inclusion of
the Perfectly Matched Layer absorbing condition~\cite{Berenger} into
the UWVF was also described.  In addition, several computational
heuristics relevant to our study were also presented.  As a result of
this work, a first parallel implementation of the UWVF called
\emph{ParMax} was written.  This was further developed at Kuava Inc,
and the University of Eastern Finland and is the basic software used
in this paper.  \emph{ParMax} uses MPI and domain decomposition to
implement parallelism.

From the point of view of theoretical convergence analysis, it was
shown in \cite{buf07} that the UWVF for the related Helmholtz equation
is a special discontinuous Galerkin (DG) method, and using duality
theory convergence estimates could be proved.  A more general DG
approach, again for the Helmholtz equation, is taken in \cite{git07}.
Based on these studies, and following the proof of explicit stability
bounds for the interior impedance boundary value problem for Maxwell's
equations in \cite{moi11}, convergence of Trefftz DG methods (in
particular the UWVF) for Maxwell's equations under rather strict
geometric assumptions was proved in~\cite{HMP13}.  A key point in this
analysis is that it is not necessary to use tetrahedral elements, but
a wide class of element shapes are theoretically covered and we shall
return to this point in Section~\ref{TISC}.

A restriction on the use of the UWVF is that the relative electric
permittivity denote $\epsr$ and relative magnetic permeability $\mur$
must be piecewise constant (constant on each tetrahedron in the
mesh). This assumption could perhaps be weakened using generalised
plane waves~\cite{IM-sylvand-21} or embedded Trefftz
techniques~\cite{PS_22} but these have not yet been extended to
Maxwell's equations (the latter technique has been demonstrated in the
open source package~\cite{NGSTrefftz}).  The restriction of piecewise
constant media is relaxed if a Perfectly Matched Layer (PML) is used.
There $\epsr$ and $\mur$ are spatially varying
tensors~\cite{Huttunen2007}.

The main issue facing the UWVF is that the condition number of the
global system rises rapidly as the number of plane wave directions
increases.  This in turn causes the iterative solution of the linear
system to slow down.  We adopt the approach of \cite{Huttunen2007}
choosing the number of plane wave directions on a given element on the
basis of its geometric size (in wavelengths) to control
ill-conditioning.

We use the same stabilised BiConjugate Gradient (BICGstab) scheme as
in \cite{Huttunen2007} for solving the global linear system resulting
from the UWVF.  Interesting new results on preconditioned iterative
schemes can be found in \cite{sirdey22}. These results are based on
the use of a structured hexahedral grid, whereas we focus on
unstructured grids in this paper.

Within the UWVF scheme used here there is considerable scope for
different choices of basis functions provided they form a complete
family of solutions of Maxwell's equations.  For example, plane wave
basis functions (as used here), vector wave functions or the method of
fundamental solutions (MFS) are discussed in \cite{badics14}.  In the
interesting case of MFS, \cite{hafner85} gives an application to wave
guides. to the method.  Our choice of plane waves in \emph{ParMax}
follows~\cite{cessenat_phd} and has the advantage that necessary
integrals of products of basis functions can be performed in closed
form on flat triangular faces in the mesh (the majority of
faces). This greatly speeds up the assembly of matrices compared to
quadrature that is needed for other choices of basis functions.

In this paper, $i=\sqrt{-1}$ and we use the convention that the time
variation of the fields and sources is proportional to $\exp(-i\omega
t)$ where $\omega$ is the angular frequency of the radiation, and $t$
is time.  All results are then reported in the frequency domain.  Bold
face quantities are vector valued. The coefficients $\epsilon_0$ and
$\mu_0$ are the permittivity and permeability of free space
respectively. The wave number $\kappa$ of the radiation is given by
$\kappa=\omega\sqrt{\epsilon_0\mu_0}$.

To further fix notation and context, we now define the Maxwell system
under consideration in this paper. Let $\Omega$ denote a Lipschitz
bounded computational domain having unit outward normal $\bfnu$ and
boundary $\Gamma:=\partial\Omega$. For a smooth enough vector field
$\bfv$, we define the tangential component $\bfv_T$ on $\Gamma$ by
$\bfv_T=(\bfnu\times\bfv)\times\bfnu$.  Then, given the wave number
$\kappa>0$, a tangential boundary vector field $\bfg$, piecewise
constant functions $\epsr$ and $\mur$, and a parameter
$Q\in\mathbb{C}$ with $|Q|\leq 1$ we seek the complex valued vector
electric field $\bfE$ that satisfies
\begin{subequations}
\begin{align}
\curl\mur^{-1}\curl \bfE-\kappa^2\epsr\bfE&=0\mbox{ in }\Omega,\label{eq:maxwell}\\
\bfnu\times\mur^{-1}\curl \bfE+\frac{i\kappa}{Z}\bfE_T&=Q\left(-\bfnu\times\mur^{-1}\curl \bfE+\frac{i\kappa}{Z}\bfE_T\right)+\bfg\mbox{ on }\Gamma.
\label{eq:bc}
\end{align}
\end{subequations}
Here $Z$ is the surface impedance (a positive real parameter).  The
boundary condition (\ref{eq:bc}) is of impedance type and is well
suited to the UWVF.  When $Q=-1$ it gives a rotated version of the
Perfectly Electrically Conducting (PEC) boundary condition for the
scatterered wave
\[
\bfE_T=\frac{Z}{2i\kappa}\bfg,
\]
where we take $g=-2i\kappa\bfE^i_T/Z$ and $\bfE^i$ is the incident
wave.  When $Q=0$ it corresponds to an outgoing condition that can be
used as a low order radiation condition, while for $Q=1$ we have a
symmetry boundary condition.

This paper presents several novel extensions of the basic UWVF that
are of considerable utility in practical applications.  In particular:
\begin{itemize}
    \item The original UWVF uses a tetrahedral grid, however the
      error estimates in \cite{HMP13} hold for more general element
      types. Besides tetrahedral elements, we have implemented
      hexahedral and wedge elements. In this paper, hexahedral
      elements are only used in the PML region.
    \item Very often curved surfaces appear in applications, and we
      have implemented a mapping technique to approximate smooth
      curved surfaces.  This allows us to use larger elements near a
      curved boundary. We shall show, using numerical experiments,
      that this improves the efficiency of the software by decreasing
      the overall time to compute a solution.  Note that we only need
      to map faces in the mesh which simplifies the implementation,
      but we have to use quadrature to evaluate integrals on curved
      faces.
    \item We show how to implement resistive sheet transmission
      conditions across thin interfaces.  Related to this we have also
      implemented a combined total field and scattered field
      formulation to allow the solution of problems involving
      penetrable media.
\item We point out that a low memory version can be used to solve very
  large problems by avoiding the storage of the most memory intensive
  matrix in the algorithm.
    \item We provide numerical results to justify the utility of the
      above innovations.
\end{itemize}

The paper proceeds as follows.  In Section \ref{sec:base} we start
with a brief derivation of the basic UWVF and describe the plane wave
based UWVF. Then in Section~\ref{TISC} we describe the five
contributions of this paper.  We start with comments on new geometric
element types and a brief discussion of implementing an algorithm for
using scattered or total fields in different subdomains in the context
of scattering from a penetrable object.  Then we move to discuss
resistive sheets, curved elements and quadrature, and finally a lower
memory version of the method.  Section ~\ref{sec:NE} is devoted to
numerical examples illustrating the UWVF and the previously mentioned
modifications.  In Section~\ref{sec:PEC} we start with two examples of
scattering from a PEC scatterer.  The first PEC example is scattering
from a sphere for which the Mie series gives an accurate solution for
comparison (c.f. \cite{Monk03}).  This example demonstrates the
benefits of curved elements and different element types, as well as
the use of very coarse meshes (compared to those used by finite
element methods). The second example is X-band scattering from an
aircraft model.  Here again we use curved elements, but in addition
use the low-memory version of the software. In Section~\ref{sec:res}
we consider two examples involving a resistive sheet.  The first is a
classical example having an exact analytical solution, and the second
is a resistive screen surrounding a sphere.  Next, in
Section~\ref{sec:het}, we consider heterogeneous or penetrable
scatterers.  The first example is a dielectric sphere where we use the
scattered/near field formulation and compare to the Mie series
solution.  The second example is scattering from a plasma.  In
Section~\ref{sec:conc} we present our conclusions. Finally in
Appendix~\ref{sec:app} we give an update to the basis selection rule
used previously in~\cite{Huttunen2007} for the new element types.
Then, in Appendix~\ref{sec:app_comp}, we present a comparison of
\emph{ParMax} with the edge finite element method for scattering from
a dielectric sphere.

\section{Derivation and properties of the basic UWVF}\label{sec:base}

In this section we provide a sketch of the derivation of the UWVF
sufficient to allow us to present the new features of this paper in
the following section.  For full details see
\cite{cessenat_phd,Huttunen2007}.

\subsection{A brief derivation of the UWVF}

The version of the UWVF presented here is equivalent to that used in
\emph{ParMax} (from \cite{cessenat_phd}) but with simplified
notation. Consider a mesh of $\Omega$ of elements of maximum diameter
$h$ denoted by ${\cal T}_h$.  An element $K\in {\cal T}_h$ in this
mesh is a curvilinear polyhedron (curvilinear tetrahedron, wedge, or
hexahedron) with boundary denoted by $\partial K$ and unit outward
normal $\bfnu^K$. We now extend the parameter $Z$ to a real piecewise
positive constant defined on all faces in the grid.  Following
\cite{cessenat_phd}, we choose $Z$ as follows. Let
\[\hat{\epsilon}=\left\{\begin{array}{rl}
|\sqrt{\epsr|_K\,\epsr|_{K'}}|& \mbox{on } K\cap K' \mbox{ for }K,K'\in{\cal T}_h\\
|{\epsr}|&\mbox{on boundary faces}\end{array}\right.
\]
where $|_K$ denotes the restriction to $K$. The edge function
$\hat{\mu}$ is defined in the same way.  Then
$Z=\sqrt{\hat{\mu}}/\sqrt{\hat{\epsilon}}$.
    
Suppose $\bfxi$ is a smooth solution of the adjoint Maxwell equation
in $K$:
\begin{equation}
\curl\overline{\mur}^{-1}\curl \bfxi-\kappa^2\overline{\epsr}\bfxi=0.
\label{max_adj}
\end{equation}
Then taking the dot product of (\ref{eq:maxwell}) with $\bfxi$ (including complex conjugation) and integrating by parts twice provides the following fundamental relation between the electric and magnetic fluxes on $\partial K$:
\begin{equation}
\int_{\partial K}{\bfnu^K\times\mur^{-1}\nabla\times\bfE}\cdot\overline{\bfxi}_T+{\bfnu^K\times\bfE}\cdot(\overline{\overline{\mur}^{-1}\nabla\times\bfxi})_T\,dA=0.\label{fund_eq}
\end{equation}
Using the above fundamental identity, we can then prove equality
(\ref{iso}) by expanding both sides of (\ref{iso}) and using
(\ref{fund_eq}). Equality (\ref{iso}) gives the conclusion of the
``isometry lemma'' (c.f. \cite{cessenat_phd}).

\begin{eqnarray}
&&\int_{\partial K}Z\left(-
\bfnu^K\times\mur^{-1}\nabla\times\bfE+\frac{i\kappa{}}{Z} \bfE_T\right)\cdot\left(-\overline{\bfnu^K\times{\overline{}}\mur^{-1}\nabla\times\bfxi+\frac{i\kappa{}}{Z} \bfxi_T}\right)\,dA.\nonumber\\
&=&\int_{\partial K}Z\left(
\bfnu^K\times\mur^{-1}\nabla\times\bfE+\frac{i\kappa{}}{Z} \bfE_T\right)\cdot\left(\overline{\bfnu^K\times{}\mur^{-1}\nabla\times\bfxi+\frac{i\kappa{}}{Z} \bfxi_T}\right)\,dA\label{iso}
\end{eqnarray}

To simplify the presentation,  we define rescaled versions of the unknowns in~\cite{cessenat_phd} as follows:
\begin{eqnarray}
\bfchi_K&=&-\bfnu^K\times\mur^{-1}\nabla\times\bfE|_{K}+\frac{i\kappa{}}{Z} (\bfE|_K)_T \label{eq:chi}\\
{\cal Y}_K&=&-\bfnu^K\times\mur^{-1}\nabla\times\bfxi|_{K}+\frac{i\kappa{}}{Z} (\bfxi|_K)_T.
\end{eqnarray}
The next step is to rewrite  (\ref{iso}) using the above quantities.  In doing so we will use the function space
of surface vectors $\bfL_T^2(\partial K):=\{\bfu\in\bfL^2(\partial K)\;|\;\bfu\cdot\bfnu^K=0\}$.
Recalling that $\bfxi|_K$ satisfies the adjoint Maxwell system in $K$ we can define $\bfF_K: \bfL_T^2(\partial K)\to  \bfL_T^2(\partial K)$ by setting
\[
\bfF_K{\cal Y}_K=\bfnu^K\times\mur^{-1}\nabla\times\bfxi|_{K}+i\kappa{}Z (\bfxi|_K)_T.
\]
Now suppose elements $K$ and $K'$ meet at a face in the mesh.  On that face $\bfnu^K=-\bfnu^{K'}$.  Also we have transmission condition requiring continuity of $\bfE_T$ and $(\mur^{-1}\nabla\times\bfE)_T$ across the face so
\begin{equation}
    \bfnu^K\times\mur^{-1}\nabla\times\bfE|_K+i\omega{}Z (\bfE|_K)_T
    =\bfchi_{K'}.\label{trans}
\end{equation}
For a boundary face we can use the boundary condition (\ref{eq:bc}) 
to replace the corresponding term on that face in terms of $\bfchi_K$ and $\bfg$.
Using the above results, we may rewrite (\ref{iso}) for every element $K$.  We conclude that for $\bfchi_K\in \bfL_T^2(\partial K)$  the following equation
holds 
\begin{align}
\int_{\partial K}Z\,\bfchi_K\cdot{\cal Y}_K\,dA=&\sum_{K'\not= K}\int_{\partial K\cap \partial K'}Z\,\bfchi_{K'}\cdot \bfF_K{\cal Y}_K\,dA\nonumber\\
+\int_{\partial K\cap \partial K'}&Z\left[Q\bfchi_K+\bfg\right]\cdot \bfF_K{\cal Y}_K\,dA,\label{eqUWVFeq}
\end{align}
for any ${\cal Y}_K\in \bfL^2_T(\partial K)$.
This equation should hold for every $K\in {\cal T}_h$, and gives the UWVF for Maxwell's equations before discretization.  Cessenat~\cite{cessenat_phd} proves the uniqueness of the UWVF solution, and existence follows because (\ref{eq:maxwell})-(\ref{eq:bc}) is well posed.

\subsection{The plane wave  UWVF (PW-UWVF)}\label{Deriv}
Now that we have the variational formulation (\ref{eqUWVFeq})
we can discretize it by using a subspace of $\bfL^2_T(\partial K)$ on each element.  It is important for efficiency that $\bfF_K$ be easy to compute 
and this is where the plane wave basis is useful~\cite{cessenat_phd}.  On each element $K$ we 
choose independent directions $\{\bfd_{K,j}\}_{j=1}^{p_K}$, $\Vert\bfd_{K,j}\Vert=1$, using the first $p_K$ Hammersley points~\cite{hardin_16} on the unit sphere.  Then we 
choose $\bfxi|_K$ to be a linear combination of the $p_K$ plane wave solutions of the adjoint problem
 \[
 \bfxi_{K,\ell,m}=\bfA_{K,\ell,m} \exp\left(i\kappa\sqrt{\overline{\epsr|_K}\overline{\mur|_K}} \bfd_{K,\ell}\cdot (\bfx-\bfx_{K,0})\right)
 \]
 for $1\leq \ell\leq P_K,$ $1\leq m\leq 2,$ and $ K\in {\cal T}_h$.  Here $\bfx_{K,0}$ is the centroid of the element and the polarizations $\bfA_{K,\ell,m}$    are chosen to be unit vectors such that $\bfd_{K,\ell}\cdot \bfA_{K,\ell,m}=0$, $m=1,2$ and $\bfA_{K,\ell,1}\cdot\bfA_{K,\ell,2}=0$.   
Now we can define a discrete subspace $\bfW_{K,h}\subset\bfL_T^2(\partial K)$     
by first defining
\begin{eqnarray*}
&&\bfW_{K,h,p_K}={\rm span}\{-
\bfnu^K\times\mur^{-1}\nabla\times\bfxi_{K,\ell,m}\\&&\qquad+(i\kappa/Z)(\bfxi_{K,\ell,m})_T, \;m=1,2,\; 1\leq \ell\leq p_K\}.
\end{eqnarray*}
Then $\bfW_{h,\bfp}=\Pi_{K\in{\cal T}_h}\bfW_{K,h,p_K}$  where 
$\bfp$ denotes the vector of number of directions on each element.

The dimension of $W_{h,\bfp}$ is $\Ndof=2\sum_{K\in{\cal T}_h}p_K$. In our work, $p_K$ is chosen according to the 
heuristic in \cite{Huttunen2007} (see Appendix \ref{sec:app} for updates to this formula for larger numbers of directions and new element types).  Then $\bfF_K$ is easy to compute using the definition
of the basis functions.  The discrete PW-UWVF uses trial and test functions from $\bfW_{h,\bfp}$ in place of $\bfW$
in (\ref{eqUWVFeq}).

In our implementation, we use domain decomposition by subdividing the
mesh according to Metis \cite{Karypis_metis98} where $p_K$ is used to estimate the work on each element.  The elements are sorted using reverse Cuthill-McKee to minimise bandwidth.  Then, enumerating the degrees of freedom element by element, the matrices and vectors corresponding
to the terms in (\ref{eqUWVFeq}) can be computed.  The left hand side of (\ref{eqUWVFeq}) gives rise to an 
$\Ndof\times \Ndof$ block diagonal Hermitian positive-definite matrix $D$, while the remaining sesquilinear 
forms on the right hand side of (\ref{eqUWVFeq}) give rise to a general sparse complex matrix $C$.  The data term involving $\bfg$ gives rise
to a corresponding vector $\vec{\bfb}$.  Denoting the vector of unknown degrees of freedom by $\vec{\bfchi}$, we 
solve the global matrix equation
\begin{equation}
(I-D^{-1}C)\vec{\bfchi}=D^{-1}\vec{\bfb}\label{dc}
\end{equation}
using BICGstab where $D^{-1}$ can be calculated rapidly element by element.  Once $\vec{\bfchi}$ is known, the solution on each element can be reconstructed for post processing.  For unbounded scattering problems, we use either the low order absorbing condition ($Q=0$) on the outer boundary or a PML as in~\cite{Huttunen2007}.  For a scattering problem the far field pattern can the be calculated using an auxiliary surface 
containing all the scatterers in its interior in the usual way~\cite{cessenat_phd}.

\section{Towards industrial scale software}\label{TISC}
In this section we discuss the extensions to the basic UWVF in this paper.
 
\subsection{New element types}
 
In \cite{HMP13} error estimates are proved for general elements that have Lipschitz boundaries, are shape 
 regular in the sense of that paper,  and are star-shaped with respect to a  ball centred at a point in 
 the element.  This allows a wide variety of elements.  We have implemented curvilinear tetrahedral, 
 wedge, and hexahedral elements.  In \emph{ParMax} we represent the boundary of each element as a 
 union of possibly curvilinear triangles.  The assembly phase is quickest if the triangles are planar
 since then quadrature can be avoided.  
 
For this paper, we rely on COMSOL Multiphysics to generate the meshes. An example of a grid using tetrahedral, wedge, and hexahedral elements is shown in Fig.~\ref{fig:pec_sphere_grid}.  Here we use wedge and hexahedral elements in the PML, and tetrahedral elements elsewhere.
 
\subsection{Scattered/total field formulation}\label{stf}
  The numerical results we shall present are all of scattering type.  The total field $\bfE$ is 
  composed
of an unknown scattered field $\bfE^s$ and a given incident field $\bfE^i$ so $\bfE=\bfE^s+\bfE^i$.  We will use 
plane wave incident fields (but point sources or other incident fields can be used) so
\begin{equation}
\bfE^i(\bfx)=\bfp\exp(i\kappa\bfd\cdot\bfx),\label{Einc}
\end{equation}
where $\bfd\in \mathbb{R}^3$ is the direction of propagation and $\Vert\bfd\Vert=1$.  The vector polarisation $\bfp\in\mathbb{C}^3$ is non-zero and satisfies $\bfd\cdot\bfp=0$.

To allow the use of a PML or other absorbing boundary condition, we need to compute using the scatterered field in the PML.  But inside a penetrable scatterer we need to compute with the total field so that there are no current sources in the scatterer.  This is standard for finite element methods, but not usual for the UWVF so we outline the process here. Suppose $\Omega$ is partitioned into two subdomains denoted $\Omega_-$ and $\Omega_+$  such that the PML (or a neighbourhood
  of the absorbing boundary if one is used) is contained in $\Omega_+$ where the scattered field is used and where $\epsr=\mur=1$.  The scatterer is contained 
  in $\Omega_-$ where possible $\epsr$ or $\mur$  are no longer unity and the total field is used.  Let $\Sigma$ denote the boundary between 
  $\Omega_-$ and $\Omega_+$ and assume that $\Sigma$ is contained in the interior of $\Omega$ and is exactly covered by faces
  of the mesh.  Then suppose that elements $K_-\subset\Omega_-$ and $K_+\subset \Omega_+$ meet at a face $F\subset\Sigma$. 

  Consider first $K_-$. Reviewing the derivation of the UWVF outlined in Section~\ref{Deriv} we see that we must rewrite the term in $\bfE$ on the right hand side
  of equation (\ref{iso}) in terms of $\bfchi_{K_-}$ and $\bfchi_{K_+}$. In particular using the facts that $\bfnu_{K_-}=-\bfnu_{K_+}$ and that on $K_+$ the scattered field $\bfE^s$ is  be approximated, so that in (\ref{eq:chi}) $\bfE|_{K^+}$ is replaced by $(\bfE^s+\bfE^i)|_{K^+}$, we obtain 
  \begin{eqnarray*}
 \left(
\bfnu^{K_-}\times\mur^{-1}\nabla\times\bfE|_{K_-}+\frac{i\kappa{}}{Z} \bfE|_{K_-,T}\right) &=&-\bfnu^{K_+}\times\mur^{-1}\nabla\times(\bfE^s+\bfE^i)|_{K_+}\\&& +\frac{i\kappa{}}{Z} (\bfE^s+\bfE^i)|_{K_+,T}=\bfchi_{K_+}+\bfg_{K_-,F}.
\end{eqnarray*}
Here we have also used the transmission condition that the tangential components of $\bfE$ and $\mur^{-1}\curl\bfE$ are continuous across $\Sigma$.  In this equation, the source function on $F$ associated with $K_-$ is
\[
\bfg_{K_-,F}=\left(\bfnu^{K_-}\times\mur^{-1}\nabla\times \bfE^i|_{F}+\frac{i\kappa{}}{Z} \bfE^i|_{F,T}\right).
\]
Carrying out the same procedure on $K_+$ remembering this element supports the scattered field gives another equation and source function
$\bfg_{K_+,F}$ again relating $\bfE$ and $\bfE^s$.  Thus a combined scattered and total field algorithm can be 
implemented by allowing for a source function $\bfg$ on the  internal
surfaces.

\subsection{Resistive sheet}
A similar procedure to that used to implement the scattered/total field UWVF can be used to derive appropriate modifications
to include resistive~\cite{jin_volakis} or conductive sheets~\cite{Senior}.  We only consider the resistive sheet.  Suppose now that a surface $\Sigma_r\subset\Omega$ is a resistive sheet, and that $\bfnu_r$ denotes a continuous normal to the sheet. This surface may be open or closed and could intersect the boundary.  We just assume that it is the union of a subset of the faces (possibly curvilinear) in the mesh. 
\begin{figure}
\begin{center}
    \includegraphics[width=0.5\textwidth]{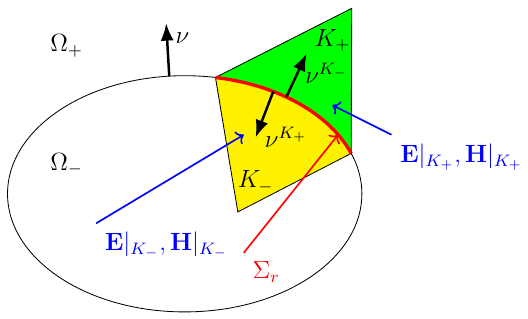} 
\caption{Geometry and notation for the resistive sheet calculation.  The normal $\nu$ is outward to $\Omega_-$.}
\label{fig:Et}
\end{center}
\end{figure}
Suppose $K_+$ and $K_-$ are two elements that meet at a face $F_r\subset \Sigma_r$ such that $\bfnu_r$ points into $K_+$.
The geometry is shown in Fig.~\ref{fig:Et}.

The resistive sheet approximation requires us to implement the following transmission conditions across $\Sigma_r$
\begin{eqnarray}
\bfnu_r\times(\mur^{-1}\curl\bfE|_{K_+}-\mur^{-1}\curl\bfE|_{K_-}) &=&i\kappa\sigma d (\bfnu_r\times \bfE_{K_+})\times \bfnu_r,\label{dit1}\\
(\bfnu_r\times \bfE|_{K_+})\times \bfnu_r&=&(\bfnu_r\times \bfE_{K_-})\times \bfnu_r,\label{dit2}
\end{eqnarray}
where $\sigma$ is the conductivity of the material in the layer and $d$ is the thickness.  Alternatively the resistivity of the sheet is given by
\[
R=(\sigma d)^{-1}=-\frac{Z_0}{i(\epsr-1)\kappa_0d}
\]
with $\epsr=1+i\sigma /(\epsilon_0\omega)$, $Z_0=\sqrt{\mu_0/\epsilon_0}$ is the impedance of free space  and $\kappa_0=\omega\sqrt{\mu_0\epsilon_0}$ is the wave-number. We now set  $\eta=\sigma d$. 

As for the scattered/total field version of UWVF we must write the term in $\bfE$ on the right hand side of (\ref{iso})
using the variables $\bfchi|_{K_+}$ and $\bfchi|_{K_-}$ on either side of $F_r$.

Adding $\bfchi|_{K_+}$ to $\bfchi|_{K_-}$, using the resistive sheet transmission conditions and noting that $\bfnu^{K_-}=\bfnu_r=-\bfnu^{K_+}$ gives
\begin{eqnarray*}
\bfchi_{K_-}+\bfchi_{K_+}
=i\kappa\left(\frac{2}{Z}\bfE_{K_+,T} + \eta \bfE_{K_+,T}\right)
\end{eqnarray*}
so
\begin{equation}
\bfE|_{K_+,T}=\bfE|_{K_-,T}=\frac{\bfchi_{K_-}+\bfchi_{K_+}}{2i\kappa/Z+i\kappa\eta }.
\label{e1nu1}
\end{equation}

Now we rewrite the term in $\bfE$ on the right hand side of (\ref{iso}) using the UWVF functions. Using Equation (\ref{e1nu1}) we have
  \begin{eqnarray}
\bfnu^{K_-}\times\mur^{-1}\nabla\times\bfE|_{K_-}+\frac{i\kappa{}}{Z} \bfE|_{K_-,T}=
\bfnu^{K_-}\times\mur^{-1}\nabla\times\bfE|_{K_-}+\frac{1}{Z} \frac{\bfchi_{K_-}+\bfchi_{K_+}}{2/Z+\eta }.
\label{eqt1}
\end{eqnarray}
But, using (\ref{dit1}), we have
\[
\bfnu^{K_-}\times\mur^{-1}\nabla\times\bfE|_{K_-}=-\bfnu^{K_+}\times\mur^{-1}\nabla\times\bfE|_{K_+}-i\kappa\eta\bfE_{K_+,T}.
\]
 Then, using the above equality in (\ref{eqt1}) together with (\ref{e1nu1}) again we have 
\begin{eqnarray*}
\bfnu^{K_-}\times\mur^{-1}\nabla\times\bfE|_{K_-}+\frac{i\kappa{}}{Z} \bfE|_{K_-,T} &=& 
-\bfnu^{K_+}\times\mur^{-1}\nabla\times\bfE|_{K_+}-i\kappa\eta\bfE_{K_+,T}+\frac{1}{Z} \frac{\bfchi_{K_-}+\bfchi_{K_+}}{2/Z+\eta }
\\
&=& 
\left(\bfchi_{K_+}-\frac{\eta}{2/Z+\eta }(\bfchi_{K_-}+\bfchi_{K_+})
\right).
\end{eqnarray*}
The new term introduces a new diagonal block into the matrix $C$ (defined before (\ref{dc})) and a perturbation to
the off diagonal blocks coupling fields on $K_+$ and $K_-$.

In the current \emph{ParMax} implementation,  $\eta$ is assumed constant on each mesh face of the resistive sheet and may be complex.

\subsection{Curved elements and quadrature}\label{CE&Q}
We take a straightforward approach to approximating curved boundaries and quadrature.
All the elements in \emph{ParMax} have faces that are unions of possibly curvilinear triangles.  Suppose a curvilinear face $F$ in the mesh is such that either an edge of $F$, or $F$ itself are entirely contained in a smooth curvilinear subset of the boundary $\Gamma$. We approximate $F$ by a mapping from a reference element $\hat{F}$ (with vertices $\widehat{\bfa}_1=(0,0)$, $\widehat{\bfa}_2=(0,1)$, and $\widehat{\bfa}_3=(1,0)$) in the $(s,t)$ plane to an approximation of $F$ using a degree
$\ell$ polynomial map $\bfF_F:\hat{F}\to F$. 

\begin{figure}
\begin{center}
    \includegraphics[width=0.6\textwidth]{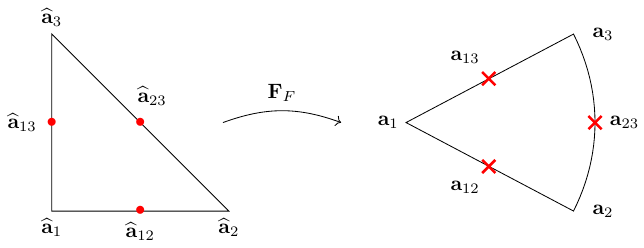} 
\end{center}
\caption{Sketch of mapping from the reference face to the face of an element in the volume mesh.  Here we sketch a quadratic map requiring that the midpoint of each edge in the curvilinear face be given.}
\label{fig:quadmap}
\end{figure}

In the important case of a quadratic map ($\ell=2$) 
we choose $\bfa_{i}$, $i=1,2,3$ to be the vertices of $K$ and take the remaining interpolation points by choosing $\bfa_{i,j}$ to be a point on the smooth boundary approximately half way between $\bfa_i$ and $\bfa_j$ (mapped from $\widehat{\bfa}_{i,j}=(\widehat{\bfa}_{i}+\widehat{\bfa}_j)/2$; see Fig.~\ref{fig:quadmap}). We use the following nodal basis
\begin{eqnarray*}
\hat{\phi}_1(s,t)&=&(1-s-t)(1-2s-2t)\\
\hat{\phi}_2(s,t)&=&s(2s-1)\\
\hat{\phi}_3(s,t)&=&t(2t-1)\\
\hat{\phi}_{1,2}(s,t)&=&4s(1-s-t)\\
\hat{\phi}_{2,3}(s,t)&=&4st\\
\hat{\phi}_{1,3}(s,t)&=&4t(1-s-t)
\end{eqnarray*}
and define $\bfF_F$ by 
\begin{equation}
\bfF_F(s,t)=\sum_{i=1}^3\bfa_{i}\hat{\phi}_i(s,t)+\sum_{i=1}^2\sum_{j=i+1}^3\bfa_{i,j}\hat{\phi}_{i,j}(s,t).
\end{equation}

We use $\bfF_F(\hat{F})\approx F$ for computing the UWVF matrices. Thus we need to evaluate integrals over each face $\bfF_F(\hat{F})$.  Using the reference element,
for a smooth function $g$ defined in a neighbourhood of $F$
\begin{eqnarray*}
\lefteqn{\int_{F}g({\bf y})dS({\bfy})\approx\int_{\bfF_F(\hat{F})}g({\bf y})dS({\bfy})}\\&=&\int_{\hat{F}} g(\bfF_F(s,t))\left\Vert \frac{\partial {\bfF_F}}{\partial s}(s,t)\times 
\frac{\partial {\bfF_F}}{\partial t}(s,t)\right\Vert\,ds\,dt.
\end{eqnarray*}

It now suffices to define quadrature on the reference element via the Duffy transform.
Suppose $G(s,t)$ is a smooth function on $\hat{F}$ (in particular the integrand above) then setting $t=(1-s)\xi$
we can map to the unit square:
\begin{eqnarray*}
\int_{\hat{F}}G\, d\hat{A}&=&\int_0^1\int_0^{1-s}G(s,t)\,dt\,ds\\&=&
\int_0^1(1-s)\int_0^1G(s,(1-s)\xi)\,d\xi\,ds.
\end{eqnarray*}
Let $(w^g_i,t_i^g)$, $1\leq i\leq N$, denote the $N$-point Gauss-Legendre rule weights and nodes on $(0,1)$. Then for each $s$
\[
\int_0^1G(s,(1-s)\xi)\,d\xi\approx \sum_{i=1}^N w_i^g G(s,(1-s)t_i^g).
\]
Next, using $N$ point Jacobi quadrature weights and nodes on $(0,1)$ denoted $(w_j^J,x_j^J)$, we obtain finally
\[
\int_{\hat{F}}G\, dA\approx\sum_{i=1}^N\sum_{j=1}^N w_j^J  w_i^g
f(x_j^J,t_i^g(1-x_j^J)).
\]
Note that the quadrature has positive weights.

\subsection{A low memory version}\label{LoMem}
A simple low memory version of \emph{ParMax} is easily available because we use BICGstab to solve the linear system.  We compute $D$ as usual (it is block diagonal),
and the vector $\vec{\bfb}$ but do not compute the elements of $C$ in (\ref{dc}).  Then, as required by BICGstab, to compute $D^{-1}C\vec\bfx$ for some vector $\vec\bfx$ we compute the blocks of $C$ element by element and accumulate $C\vec\bfx$ 
element by element. Then $D^{-1}$ is computed element by element using precomputed LU decompositions.   Obviously computing the entries of $C$ repeatedly greatly increases CPU time but this allows us to compute solutions to 
problem that would otherwise require very large memory to store $C$.  For example, the solution of a scattering problem for a full aircraft at X-band frequencies is shown in Section~\ref{Aircraft}.

\section{Numerical examples}\label{sec:NE}

All results were generated using the computer clusters Puhti and Mahti at the CSC – IT Center for Science Ltd, Finland. Detailed descriptions of these supercomputers can be found from the CSC's website \cite{CSC}. Computational grids used in this work were prepared using COMSOL Multiphysics on a personal computer. In addition, the geometry model for aircraft used in Section~\ref{Aircraft} is adapted from the COMSOL's application {\it Simulating Antenna Crosstalk on an Airplane's Fuselage}.

For all numerical experiments, the incident electric field is a plane wave propagating in the direction of the positive $x$-axis polarised in the $y$-direction where the field is given by (\ref{Einc}). 

\subsection{Scattering from PEC objects}\label{sec:PEC}

In both PEC examples, we compute the scattered field and the incident field is used as a source via the PEC condition on the surface of the scatterer.

\subsubsection{A sphere} \label{scat_sphere}

This experiment is intended to demonstrate the advantages of multiple
element types, and a cuvilinear approximation to a smooth curved
boundary (the surface of the sphere). In particular, we study
scattering from a PEC sphere placed in vacuum, $\epsr=\mur=1$, where
the frequency of the incident field is $f = 2$ GHz (wavelength in air:
$\lambda_0=0.19946$ m). The scatterer is a PEC sphere with radius of 1
meter or approximately $5\lambda_0$.

In order to demonstrate the use of several types of large elements and
curvilinear grids, we use an unusually large computational domain. In
particular, the sphere is placed with its origin at the center of a
cube shaped computational domain $[-1-15\lambda_0, 1+15\lambda_0]^3$.

An absorbing boundary condition, (\ref{eq:bc}) with $Q=0$, is used on
the exterior surface geometry. In addition, a PML with thickness of
$5\lambda_0$, and constant absorption parameter $\sigma_0=1$ is used
within each side of the cube.

\begin{figure}[!htb]
    \centering
    \includegraphics[width=0.6\textwidth]{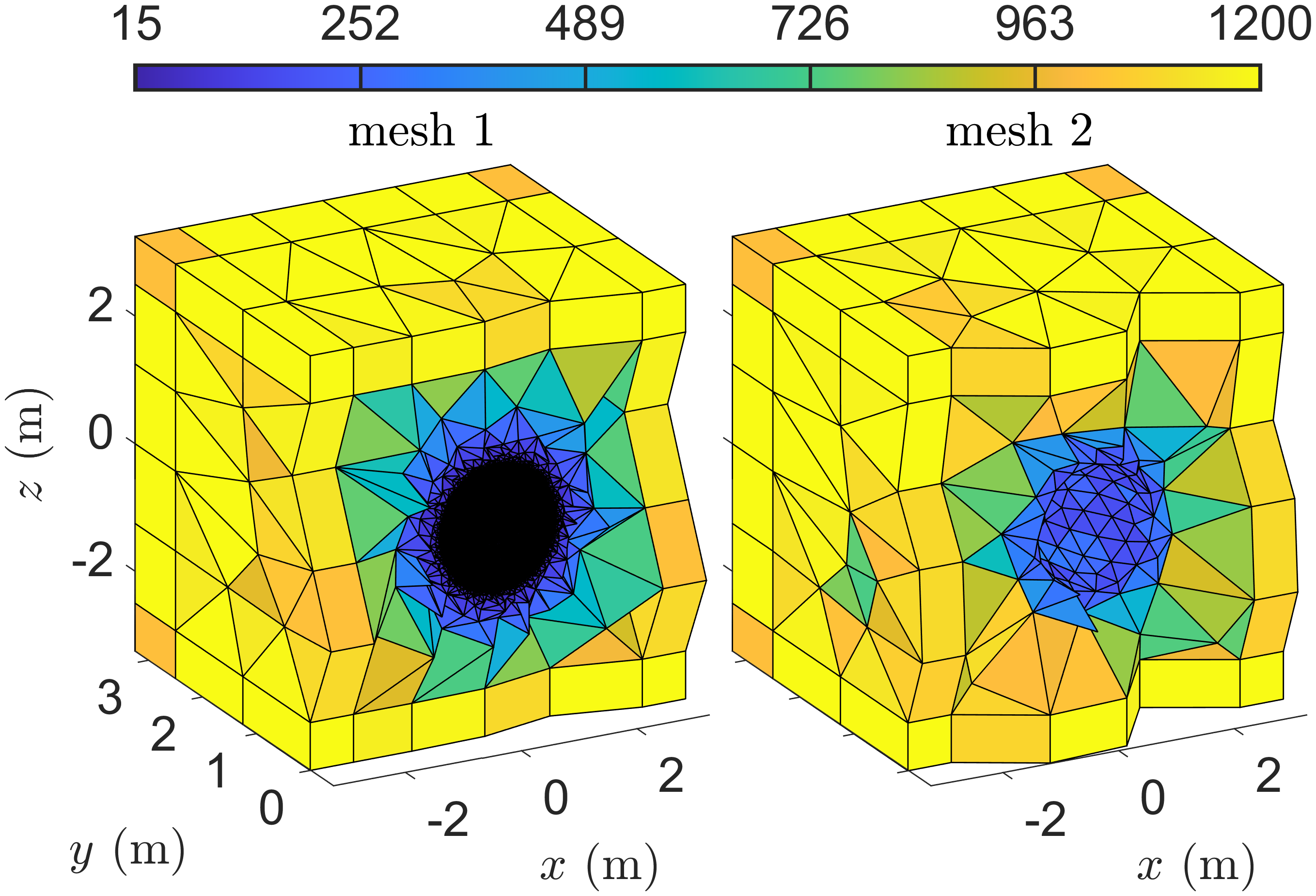} 
\caption{Cross-sections of the computational grids used to approximate scattering from a PEC sphere. The colorbar shows the number of plane wave  directions on each element.  The major difference is the grid density on the surface of the 
sphere.  This figure shows how wedge and hexahedral elements can be usefully employed in  the outer PML layer. In all figures showing grids, the colorbar shows the number of plane wave direction, see Eq. (\ref{eq:Nell}), for each element. }
    \label{fig:pec_sphere_grid}
\end{figure}

Two computational grids with different geometric approximations are
used (see Fig.~\ref{fig:pec_sphere_grid}). For the first grid
(\emph{mesh 1}), we set the requested element size on the PEC sphere
to be $\hs=\lambda_0/5$ which we will see provides a geometrically
accurate surface representation using flat face elements such that the
far field pattern predicted by UWVF is in good agreement with the far
field computed via Mie scattering. For the second grid (\emph{mesh
2}), the surface grid density is relaxed to $\hs=3\lambda_0$. For both
cases, the requested element size in the volume is set to
$\hv=10\lambda_0$. We request wedge and hexahedral elements in the PML
region, and tetrahedra elsewhere.

We show results for three cases: 1) scattering calculated using
\emph{mesh 1}, 2) scattering computed using \emph{mesh 2} with flat
facets approximating the surface of the sphere, and 3) the use of
\emph{mesh 2} with a quadratic approximation to the boundary of the
sphere (see Section~\ref{CE&Q}). These results are computed on the
Puhti system with 5 computing nodes and 10 cores per
node. Table~\ref{not} summarises the notation used in reporting
results and Table~\ref{tab:pec} gives details of the PEC computations
including CPU time.

Figure~\ref{fig:pec_sphere_grid} shows that up to 1,200 directions are
used on some elements.  As is done in the Mie series, the field
scattered by a sphere can be approximated well by relatively few
spherical vector wave functions, and these in turn can be approximated
by special plane wave expansions involving many less
directions~\cite{macphie03}.  However our code is for general wave
propagation problems, and the heuristic used in the appendix for
giving the number of directions on a element always chooses the
largest number consistent with a chosen condition number.  This
approach is intended to ensure good accuracy (within the conditioning
constraint) for a general problem.

\begin{table}[!htb]
    \caption{Notation used in reporting computational results.}
    \centering
    \begin{tabular}{r|l}
    Symbol & Definition\\\hline
        $\Nt$ & Number of tetrahedral elements\\
        $\Nw$ & Number of wedge elements\\
        $\Nh$ & Number of hexahedral elements\\
        $\Nv$ & Number of vertices\\
        $h_{\min}$  & Minimum distance between vertices\\
        $h_{\max}$  & Maximum distance between vertices\\
        $\Ndof$ & Number of degrees of freedom\\
        $N_{{\mathrm{iter}}}^{{\mathrm{BiCG}}}$ &Number of BiCGstab iterations required to reach the \\&requested tolerance value $10^{-5}$\\
        L2-error & Relative L2-error computed from the bistatic RCS \\
        CPU-time & Elapsed wall-clock time needed to \\&assemble the matrices and reach the solution 
    \end{tabular}
    \label{not}
\end{table}

\begin{figure}[!htb]
    \centering
    \includegraphics[width=0.55\textwidth]{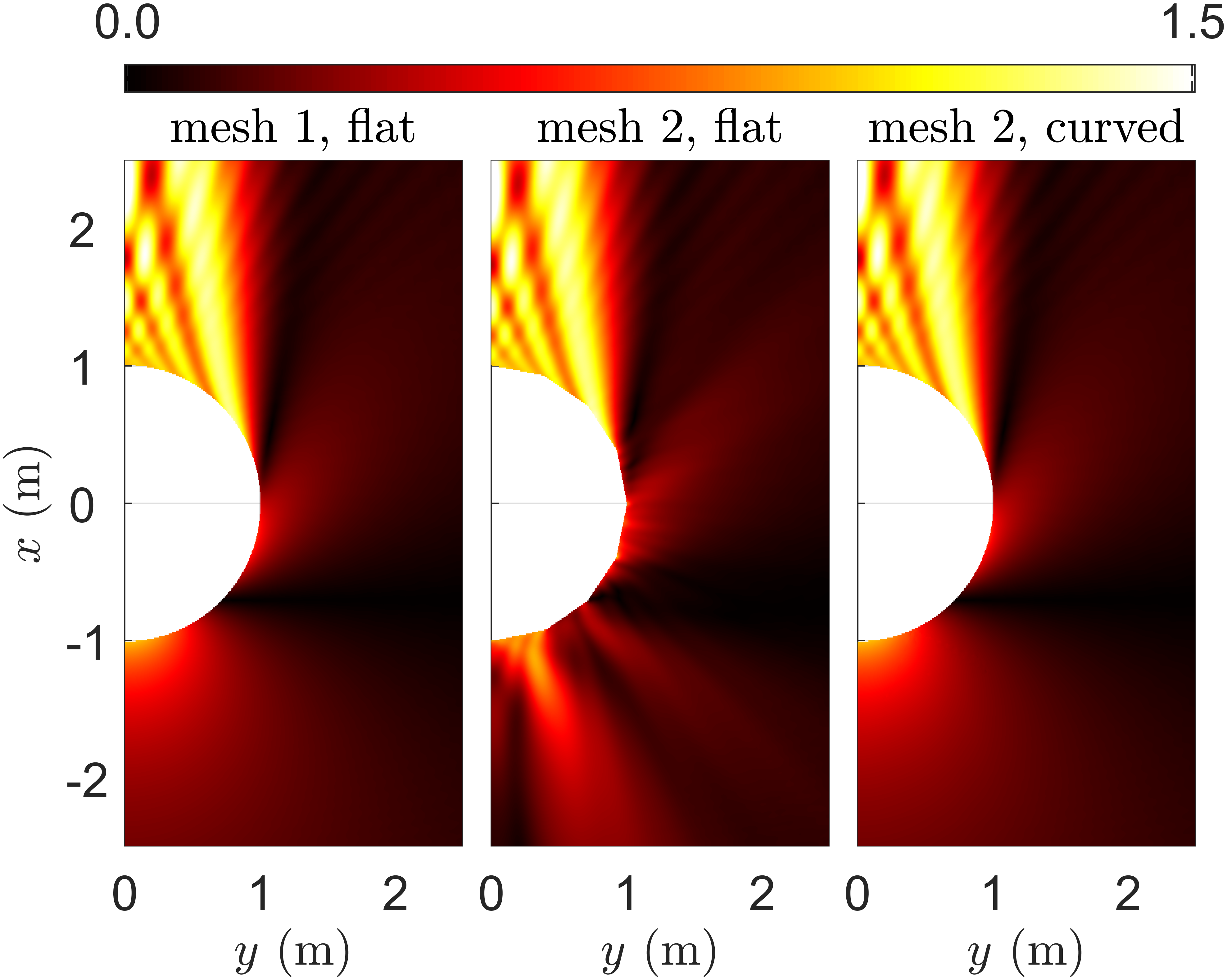} 
\caption{Snapshots of the scattered electric field component $|E_y^s|$ for the meshes considered here.  There is good agreement between \emph{mesh 1} and \emph{mesh 2} with curved faces.  \emph{Mesh 2} with flat faces produces unacceptable error due to the coarse boundary approximation.}
    \label{fig:pec_sphere_snap}
\end{figure}

Figure \ref{fig:pec_sphere_snap} shows the modulus of the $y$-component of the scattered electric field $|E_y^s|$ on the $z=0$ plane.  There are clear differences between the results for \emph{mesh 1} and \emph{mesh 2} with flat facets.  These are caused by the coarse surface grid in the second case.  However \emph{mesh 1} with flat facets, and \emph{mesh 2} with a quadratic boundary approximation are in good agreement. As can be seen in Table~\ref{tab:pec}, \emph{mesh 2} with quadratic boundary approximation is much cheaper in terms of CPU-time than \emph{mesh 1}.

Very often the far field pattern of the scattered wave~\cite{Monk03} is the quantity of interest for these calculations, and in particular the Radar Cross Section 
(RCS) derived from the far field pattern. In this paper, far field directions are defined in terms of the azimuth 
angle $\phi$ (${}^\circ$) as $(\cos(\phi\pi/180), \sin(\phi\pi/180),0)$.

\begin{figure}[!htb]
    \centering
    \includegraphics[width=0.6\textwidth]{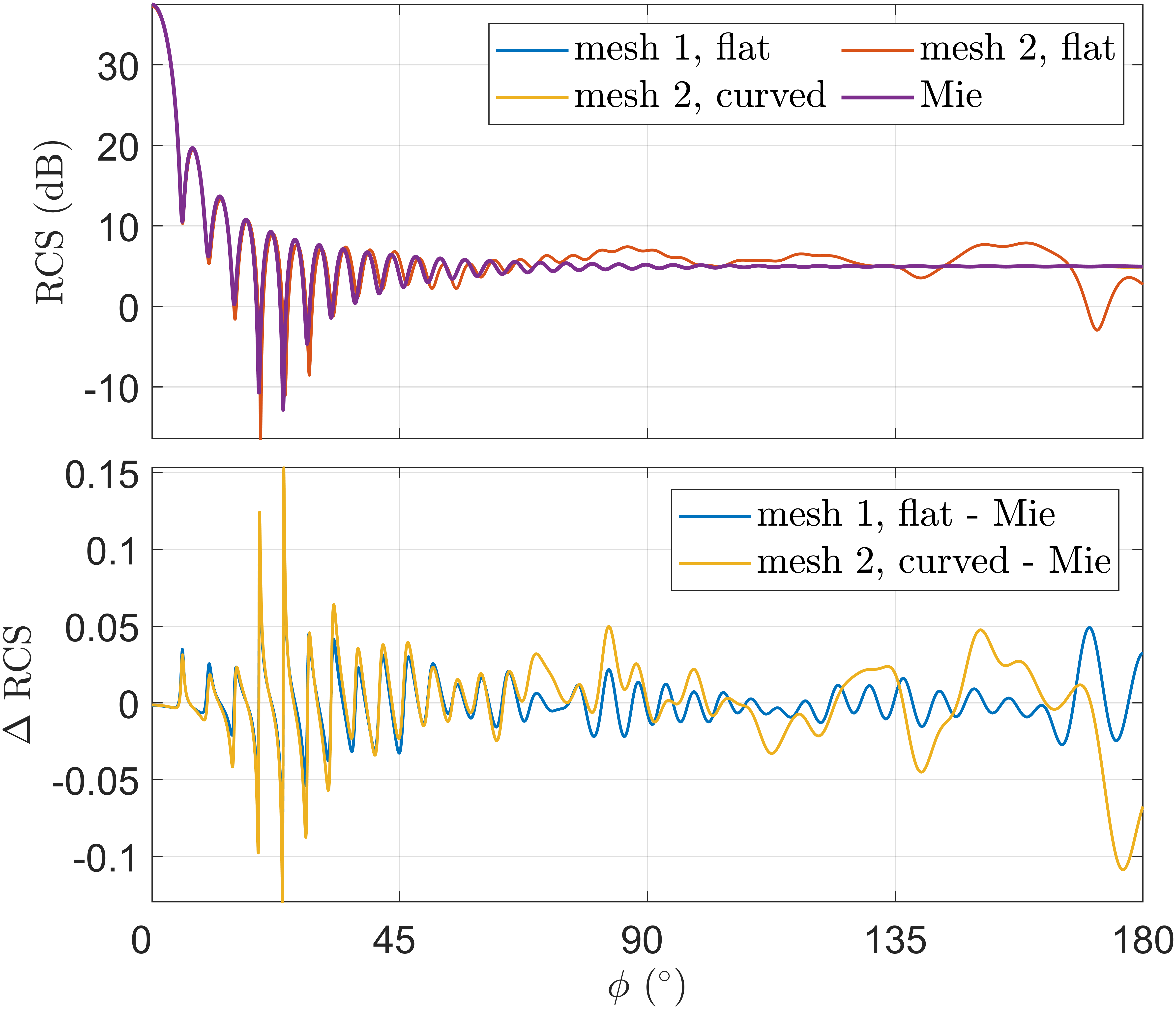} 
\caption{Bistatic RCS at 2 GHz for the PEC sphere.  We show the RCS
  computed using meshes 1 and 2 compared to the Mie series. The bottom
  panel shows the difference between the numerical solution and the
  Mie series. Note that the 'Mesh 2, flat - Mie' curve is omitted from
  the bottom panel due to its significantly larger amplitude compared
  to the other two curves.}
    \label{fig:pec_sphere_mie}
\end{figure}

Figure \ref{fig:pec_sphere_mie} shows a comparison of the bistatic RCS predicted by the computational experiments with the UWVF and one computed by the Mie series.  Clearly \emph{mesh 2} with flat facets produces an inaccurate far field pattern, whereas \emph{mesh 1} or \emph{mesh 2} with curved elements produces  much more accurate predictions.
 
\begin{table}[!htb]
    \caption{Table giving  details for the computational grids for the PEC sphere, and accuracy and timing. See Table~\ref{not} for definitions of the reported quantities. }
    \centering
    \resizebox{\textwidth}{!}{
      \begin{tabular}{cc|ccccccc|ccc}
        mesh id & surface & $\Nt$ & $\Nw$ & $\Nh$ & $\Nv$ & $h_{\min}$ (cm) & $h_{\max}$ (m) & $\Ndof$ & $N_{{\mathrm{iter}}}^{{\mathrm{BiCG}}}$ & L2-error (\%) & CPU-time (s)\\
        \hline
1 & flat & 122,680 & 228 & 56 & 29,890 & 0.65 & 1.93 & 9,325,028 & 211 & 0.21 & 730\\
2 & flat & 1,384 & 228 & 56 & 593 & 12.79 & 2.13 & 1,996,308 & 73 & 25.12 & 355\\
2 & curved & 1,384 & 228 & 56 & 593 & 12.79 & 2.13 & 1,996,308 & 72 & 0.36 & 490
\end{tabular}
}
      \label{tab:pec}
\end{table}

\subsubsection{An aircraft at X-band frequency}\label{Aircraft}

The aircraft model used in this section is derived from a model available in COMSOL (application \emph{Simulating Antenna Crosstalk on an Airplane’s Fuselage}).  We treat the aircraft as a curvilinear perfect conductor.  The frequency of the incident field is $f = 8$ GHz so $\lambda_0=0.03737$ m. The aircraft is 20.5 m or 547 wavelengths long and 17.8108 m or 475 wavelengths wide. A perfectly matched layer with a thickness of $5\lambda_0$ is added to each side of the cuboid computational domain with side lengths (16.3334, 16.3334, 4.4826) m. 

For generating the computational grid, curved elements and a mesh size
parameter of $\hs=3\lambda_0$ were employed on the aircraft's surface.
Because we can use 10$\lambda_0$ sized elements away from the
boundary, the entire grid can be created on a standard office
computer. The grid consists of 697,783 tetrahedral elements with
142,731 vertices covering the computational domain. The surrounding
PML layer is discretized using 413 hexahedral and 14,904 wedge
elements. In addition, for this grid $h_{\min} = 1.72$ cm and
$h_{\max} = 0.55$ m. The number of degrees of freedom $\Ndof$ is
  685,245,422.

\begin{figure}[!htb]
    \centering
    \includegraphics[width=0.65\textwidth]{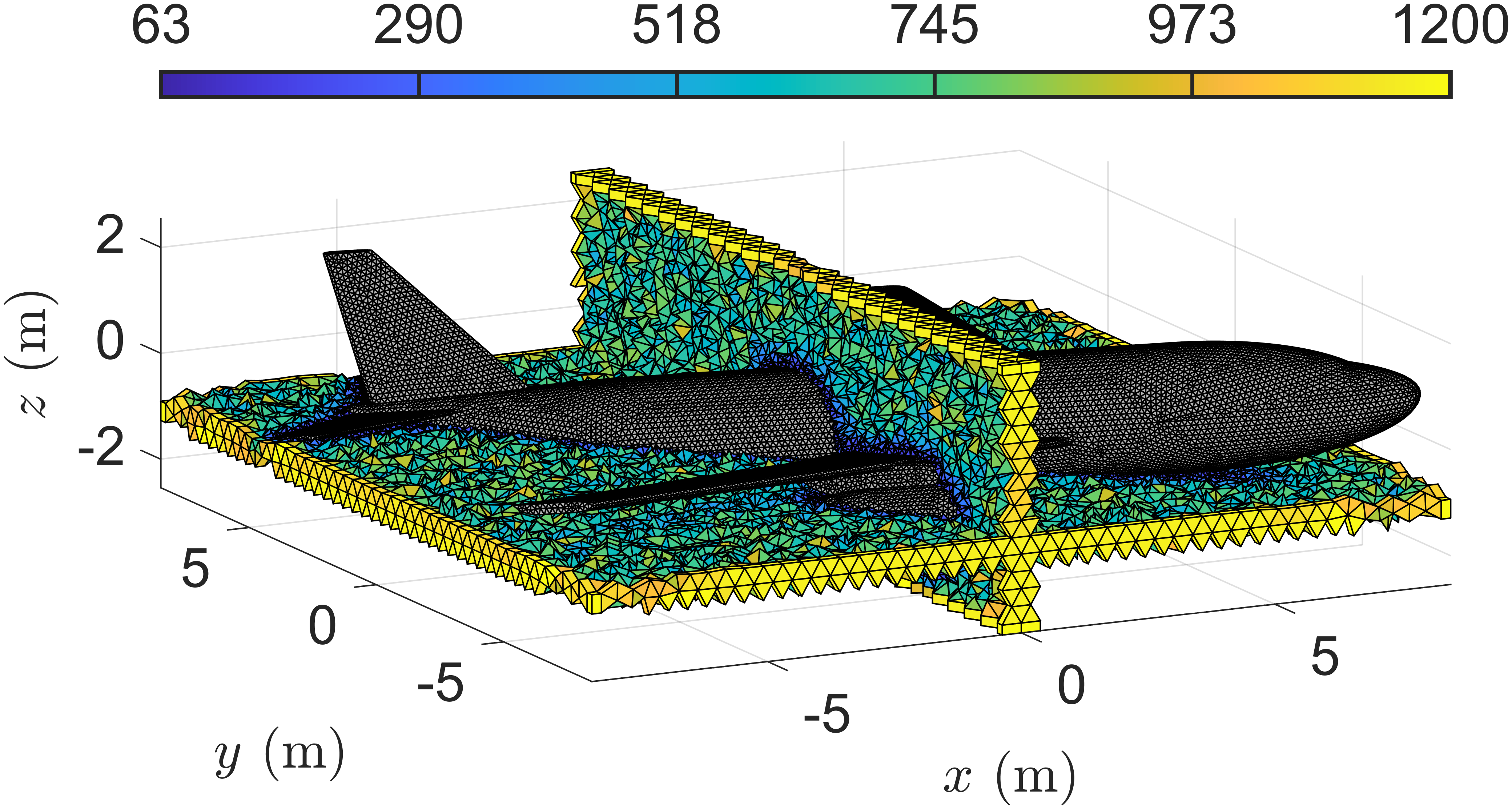} 
\caption{Surface triangulation for the aircraft and tetrahedral, hexahedral, and wedge elements on two planes. The colorbar shows the number of plane waves per element. }
    \label{fig:aircraft_grid}
\end{figure}

In Fig. \ref{fig:aircraft_grid}, the computational grid on the
aircraft's surface is depicted, along with the grid in two planes:
$z=-1$ m and $x=0$ m.

For this numerical experiment we used the supercomputer Mahti. When we
tried this example storing the matrices $C$ and $D$ as usual (see
Eq.~(\ref{dc})), we ran out of memory. An estimate for the memory
requirement for solving the problem by keeping all necessary matrices
in memory is 68 terabytes. The low memory version relaxes this memory
requirement by a factor of 0.2. So we switched to the low memory
version described in Section~\ref{LoMem} to compute the results shown
here. In Mahti, we used a total of 200 computing nodes and 100 CPU
units per node. The total time for the calculation was 18 hours,
including building the system matrix $D$, and then iteratively
reaching the requested solution accuracy. The BICGstab algorithm took
a total 437 iterations.

\begin{figure}[!htb]
    \centering
\includegraphics[width=0.5\textwidth]{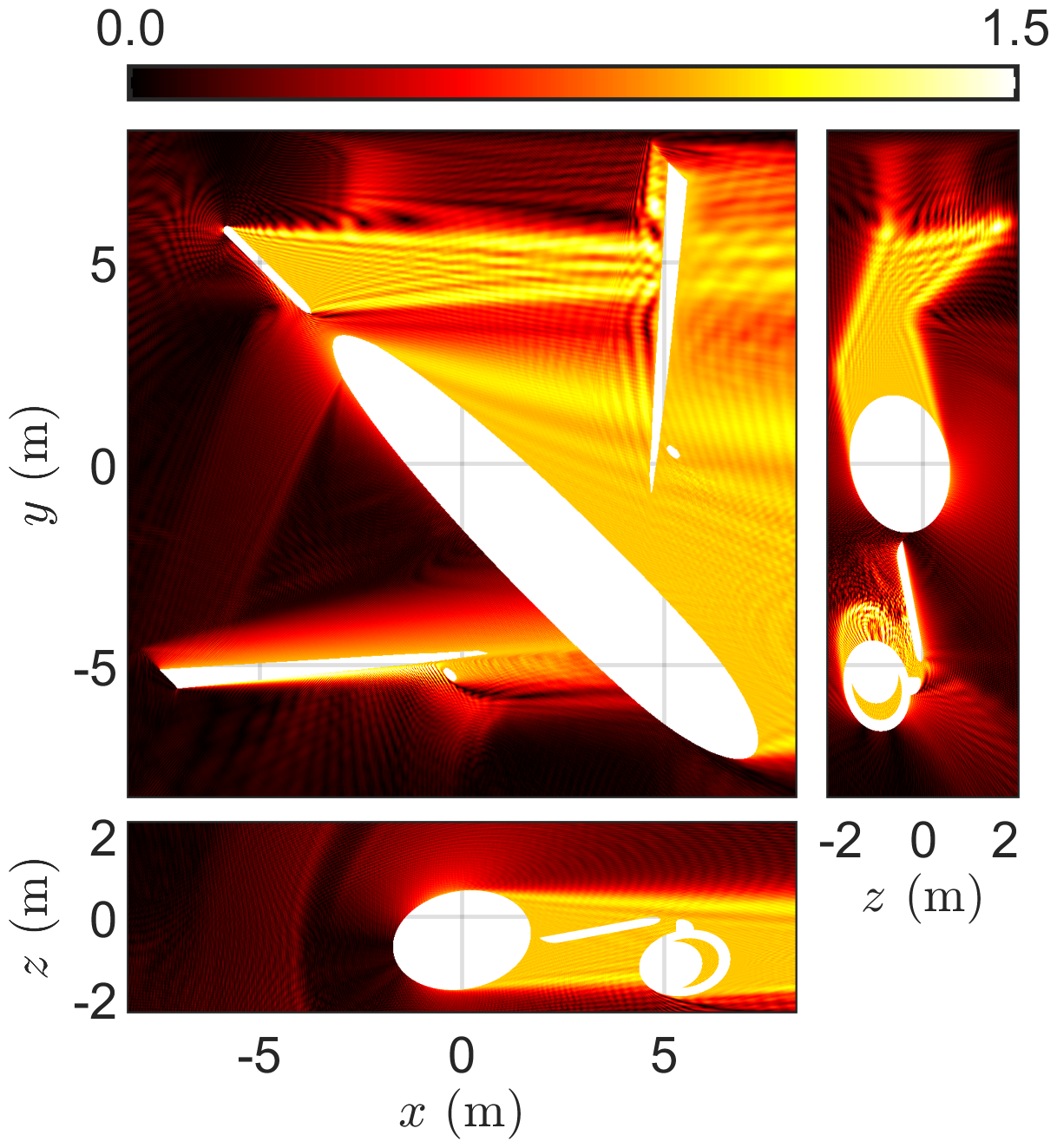} 
\caption{Snapshots of the scattered electric field $|E_y^s|$ for the aircraft model. In the top left panel we show results in the  $x-y$ plane at $z=0$, top right is in the $z-y$ plane at $x=0$ and bottom left is in the $x-y$ plane at $z=0$.  Clearly  a strong shadow region and multiple reflections are evident.}
    \label{fig:aircraft_snap}
\end{figure}

\begin{figure}[!htb]
    \centering
    \includegraphics[width=0.9\textwidth]{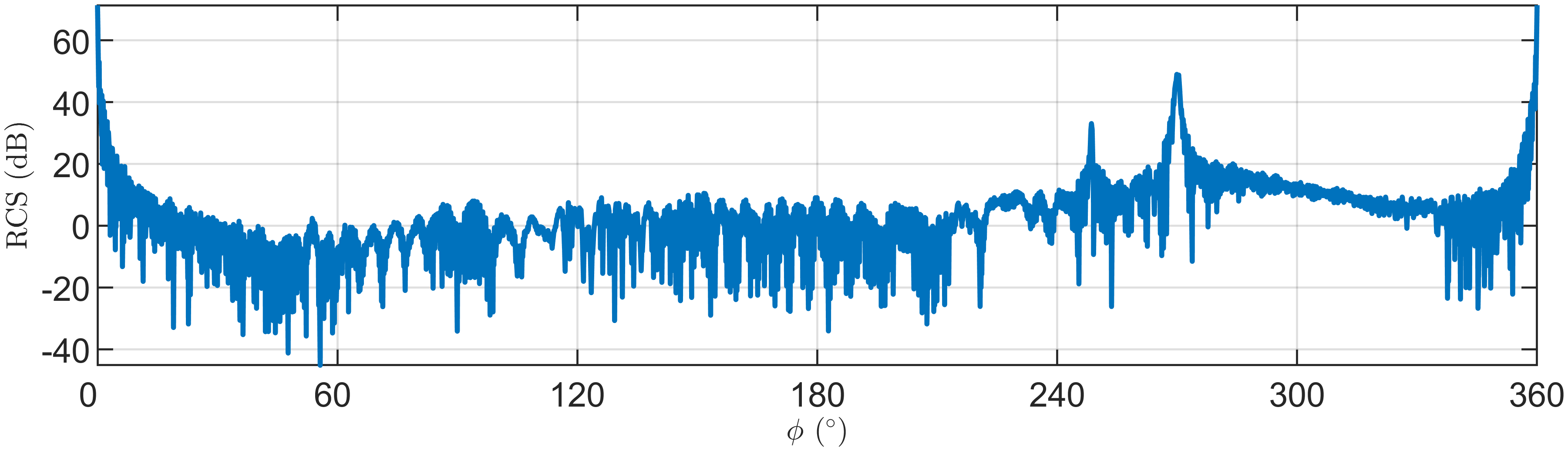}
\caption{Bistatic RCS for the aircraft at 8 GHz.}
    \label{fig:aircraft_far}
\end{figure}

Figure \ref{fig:aircraft_snap} shows the scattered electric field $|E_y^s|$ on the $x=0$, $y=0$, $z=0$ planes. Fig. \ref{fig:aircraft_far} shows the RCS in a full azimuth angle range $\phi\in[0, 360]^\circ$. 

\subsection{Resistive sheets}\label{sec:res}

\subsubsection{Salisbury screen}

We next model a ‘`Salisbury screen'’ (W. W. Salisbury, U. S. Patent US2599944 A 1952). For simplicity, assume $\epsr=\mur=1$.
Our standard incident field (\ref{Einc}) propagates normally to a resistive sheet
at $x=-H$, $H>0$ backed by a PEC surface at $x=0$.

To the left of the resistive
sheet the total electric field is 
\[
{\bf E}_{-}=\left(\begin{array}{c}0\\1\\0\end{array}\right)
\exp(i \kappa x)+ \left(\begin{array}{c}0\\R_2\\0\end{array}\right)\exp(-i\kappa x),\; x<-H,
\]
where $(0,R_2,0)^T$ is the polarization of the reflected wave.  Between the sheet and the PEC surface, where $-H<x<0$,
\[
{\bf E}_+=\left(\begin{array}{c}0\\q_{0,2}\\0\end{array}\right)
\exp(i \kappa x)+ \left(\begin{array}{c}0\\q_{1,2}\\0\end{array}\right)\exp(-i\kappa x).
\]
Here $(0,q_{0,2},0)^T$ and $(0,q_{1,2},0)^T$ are the polarizations of the left and right going waves respectively in the gap $-H<x<0$.

Imposing the PEC boundary condition at $x=0$ and the resistive sheet
transmission conditions (\ref{dit1})-(\ref{dit2}) at $x=-H$ shows that
\begin{eqnarray*}
  \mathit{R_2} &=& 
-\frac{\left(i\left( \eta -1\right) \sin \! \left(\kappa H \right)-\cos \! \left(\kappa H \right)\right) \exp\left(-2 i \kappa H\right)}{i\left( \eta +1\right) \sin \! \left(\kappa H \right)-\cos \! \left(\kappa  H \right)}, \\\mathit{q_{0,2}} &=& 
\frac{\exp\left(-i \kappa H\right)}{-i(\eta+1) \sin \! \left(\kappa H \right) +\cos \! \left(\kappa H \right)},\\
 \mathit{q_{1,2}} &= &
\frac{\exp\left(-i \kappa H\right)}{i\left( \eta +1\right) \sin \! \left(\kappa H \right)-\cos \! \left(\kappa H \right)}.
\end{eqnarray*}
For given $\kappa$, $H$, and $\eta$ we can compute total field in each region and compare to the analytic solution.  As is well known, one choice of $\eta$ gives $R_2=0$:
\[
\eta=1-i\cot(\kappa H).
\]
A particularly interesting case occurs when
$\cot(\kappa H)=0$ or $H=\pi/(2\kappa) $.  Recalling that the wavelength of the radiation is $\lambda_0=2\pi/\kappa$, we see that zero reflection occurs when $H=\lambda_0/4$.

\begin{figure}[!htb]
    \centering
    \includegraphics[width=0.55\textwidth]{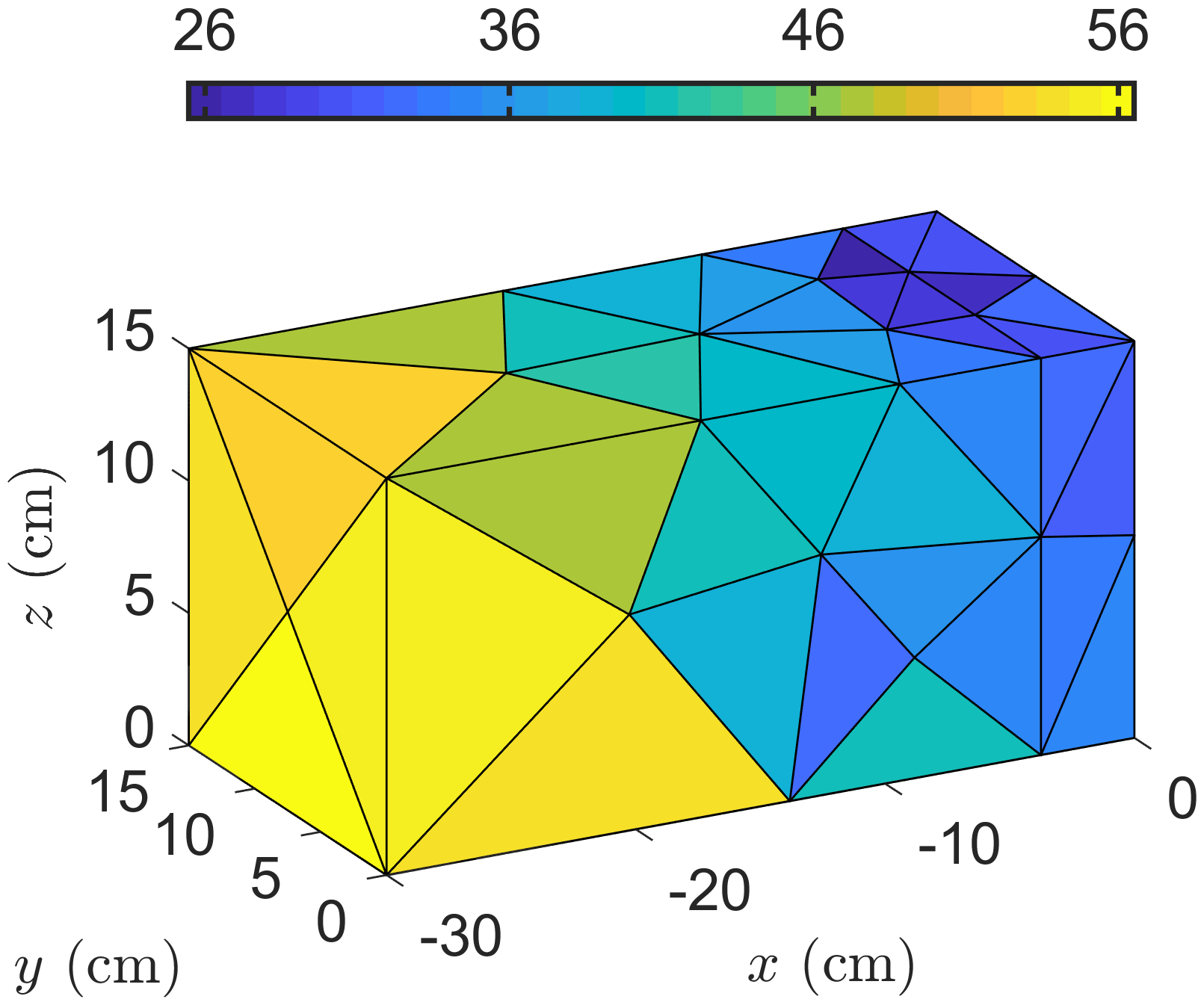} 
\caption{Computational grid for the Salisbury screen, showing smaller elements between the screen and PEC surface, expanding away to the left.}
    \label{fig:sheetpw_grid}
\end{figure}

To test the UWVF for resistive sheets,  we use a rectangular parallelepiped computational region with faces normal to the coordinate directions (see Fig.~\ref{fig:sheetpw_grid}).  The rightmost face $x=0$ is PEC ($Q=-1$) the left-most face is an ABC with $Q=0$ where we use a non-homogeneous absorbing boundary condition to excite the plane wave. There is no need for a PML since the solution is a wave propagating orthogonally to the absorbing boundary.

 On the remaining faces $Q=\pm 1$ chosen so that the incident plane wave propagates along the box without distortion. We take the radiation to have frequency $f=2$ GHz and place the  resistive sheet one quarter wavelength ($H=3.75$ cm) from the PEC surface.  Figure \ref{fig:sheetpw_grid} shows a surface of the grid used in the computations. In this case, the grid consists of 136 tetrahedral elements with 51 vertices. In addition, $h_{\min} = 3.33$ cm, $h_{\max} =  16.18$ cm, and $\Ndof=9,946$.

\begin{figure}[!htb]
    \centering
\includegraphics[width=0.55\textwidth]{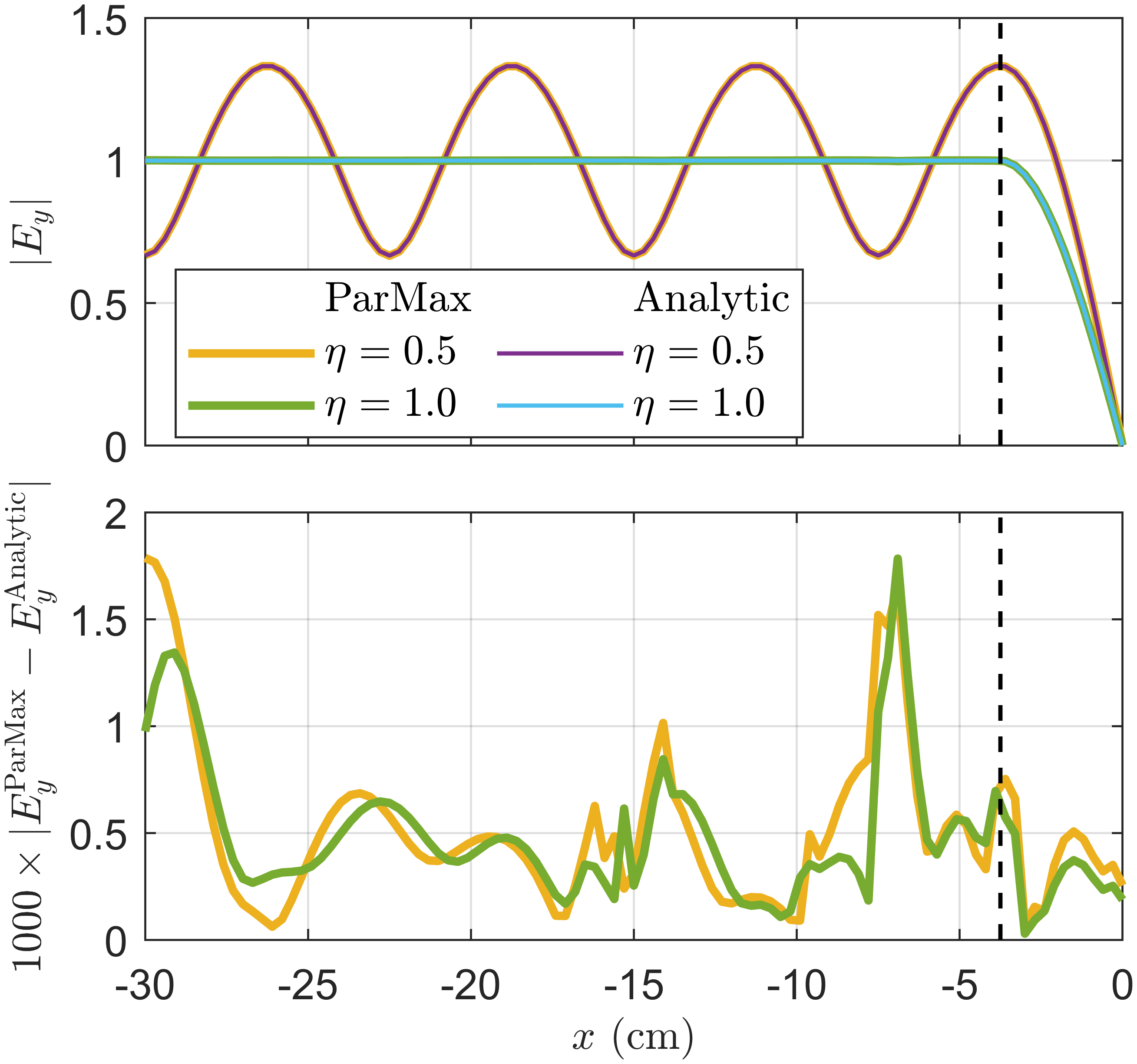} 
\caption{Magnitude of the $y$-component of the total electric field
  $E_y$ as a function $x$ with $y=z=7.5$ cm with $\eta=1$ and
  $\eta=0.5$, together with the analytic solution.  The vertical
  dashed line marks the location of the resistive sheet.  To the left
  of the sheet, no reflected wave is evident when $\eta=1$ whereas a
  reflected wave is indicated when $\eta=0.5$. The bottom panel shows
  the absolute value of the difference between the numerical and
  analytic solutions, with the difference multiplied by a factor of
  1,000.}
    \label{fig:sheetpw_snap}
\end{figure}

Results are shown in Fig. \ref{fig:sheetpw_snap}.  We plot the magnitude of the $y$-component of the total field as a function of $x$ when  $y=z=7.5$ cm.  Choosing $\eta=1$ the magnitude of the field is flat to the left of the resistive sheet  showing that this choice of $\eta$ gives rise to no reflected wave.  However when $\eta=0.5$ the non-constant magnitude of the total field indicates that a reflected wave is present to the left of the sheet. 

In this example just two plane waves per element would suffice to
compute the solution, but our code is not tuned to this example and
chooses the number of directions as if this is a general Maxwell
problem via the condition number heuristic in
Appendix~\ref{sec:app}. 

\subsubsection{Sphere with resistive sheet}

In this second experiment with resistive sheets, a PEC sphere with a radius of 1 meter is placed at the origin of a cube $[-1-15\lambda_0, 1+15\lambda_0]^3$. Surrounding this sphere is a spherical resistive sheet of radius $1+\lambda_0/4$ and surrounding both is an artificial sphere of radius $1+\lambda_0$. Outside this artificial sphere we  compute the scattered field, and inside the total field (see Section~\ref{stf}). The incident field then gives rise to a source  on the artificial boundary as detailed in Section~\ref{stf}. A PML with a thickness of $5\lambda_0$ is applied to the inside each side of the cube, and the frequency of the incident field is set at $f = 2$ GHz. The incident field gives rise to a source located on the sphere.

\begin{figure}[!htb]
    \centering
    \includegraphics[width=0.6\textwidth]{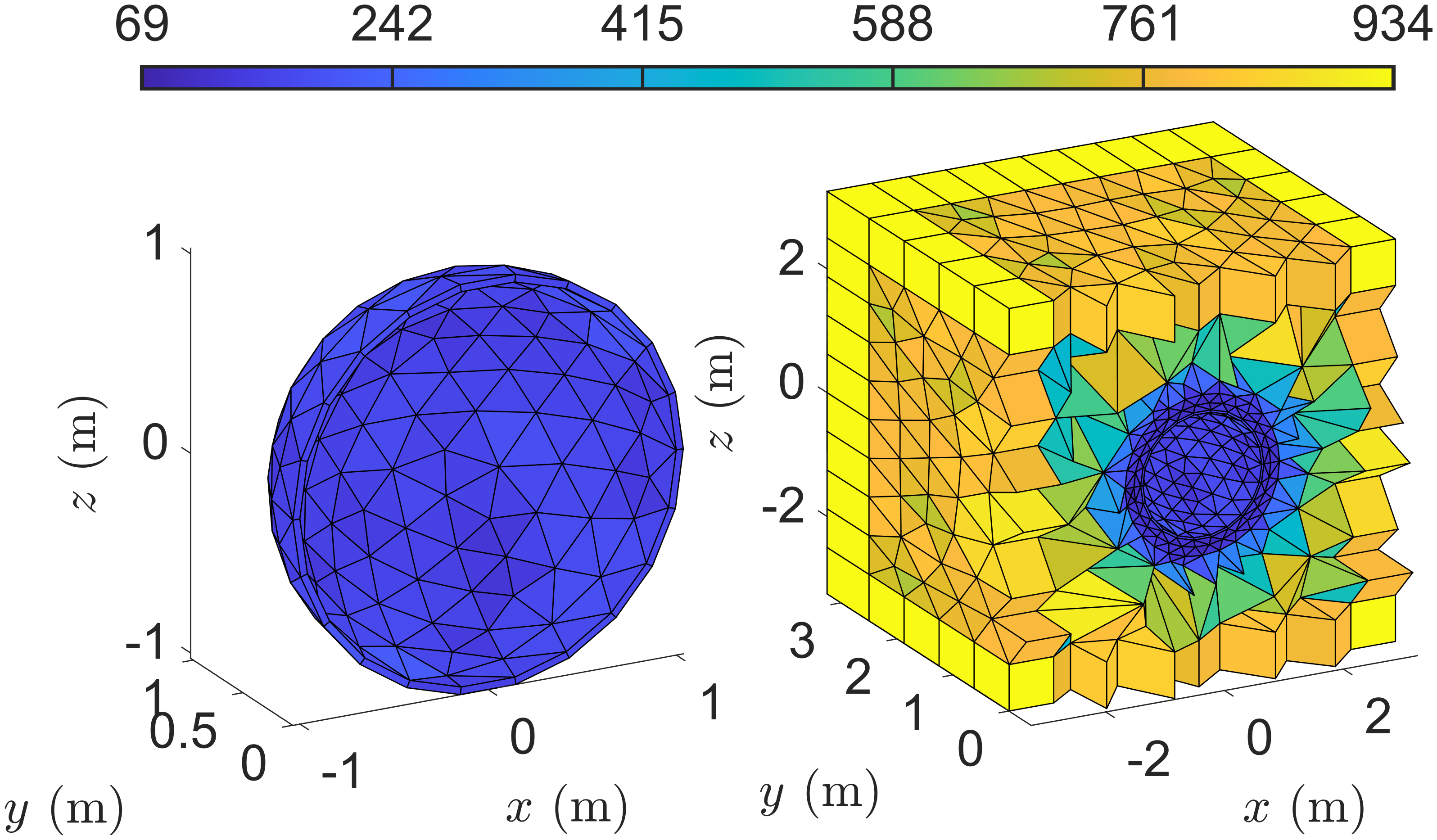} 
\caption{Left: 428 wedge elements forming the interior of the resistive sheet outside the PEC surface. Right: A cross-section of the overall computational grid. }
    \label{fig:sheet_sphere_grid}
\end{figure}

To generate the computational grid, we utilized curved elements and a
mesh size parameter of $\hs=2\lambda_0$ on the PEC and resistive sheet
surfaces. The grid comprises of 428 wedge elements (representing the
domain between the resistive sheet and PEC sphere) and 7,589
tetrahedral elements that cover the main domain of
interest. Furthermore, the surrounding PML layer is discretized using
a combination of 104 hexahedral and 936 wedge elements. The entire
grid is composed of 2,550 vertices, with $h_{\min} = 3.75$ cm and
$h_{\max} = 1.48$ m. A cross-section of the computational grid is
shown in Fig. \ref{fig:sheet_sphere_grid}. In this case, the number of
degrees of freedom $\Ndof=4,599,234$.

We used the supercomputer Puhti with a total of 5 computing nodes and 10 CPU units per node to solve the three configurations for different resistive sheet parameters $\eta$. It took 18 ($\eta=0$), 17 ($\eta=0.5$), and 17 ($\eta=1.0$) minutes CPU-time respectively to build the 
system matrices and then reach the solution with BICGstab. The solution was achieved after 194 iterations for $\eta=0$, 169 iterations for $\eta=0.5$, and 172 iterations for $\eta=1.0$.

\begin{figure}[!htb]
    \centering
    \includegraphics[width=0.55\textwidth]{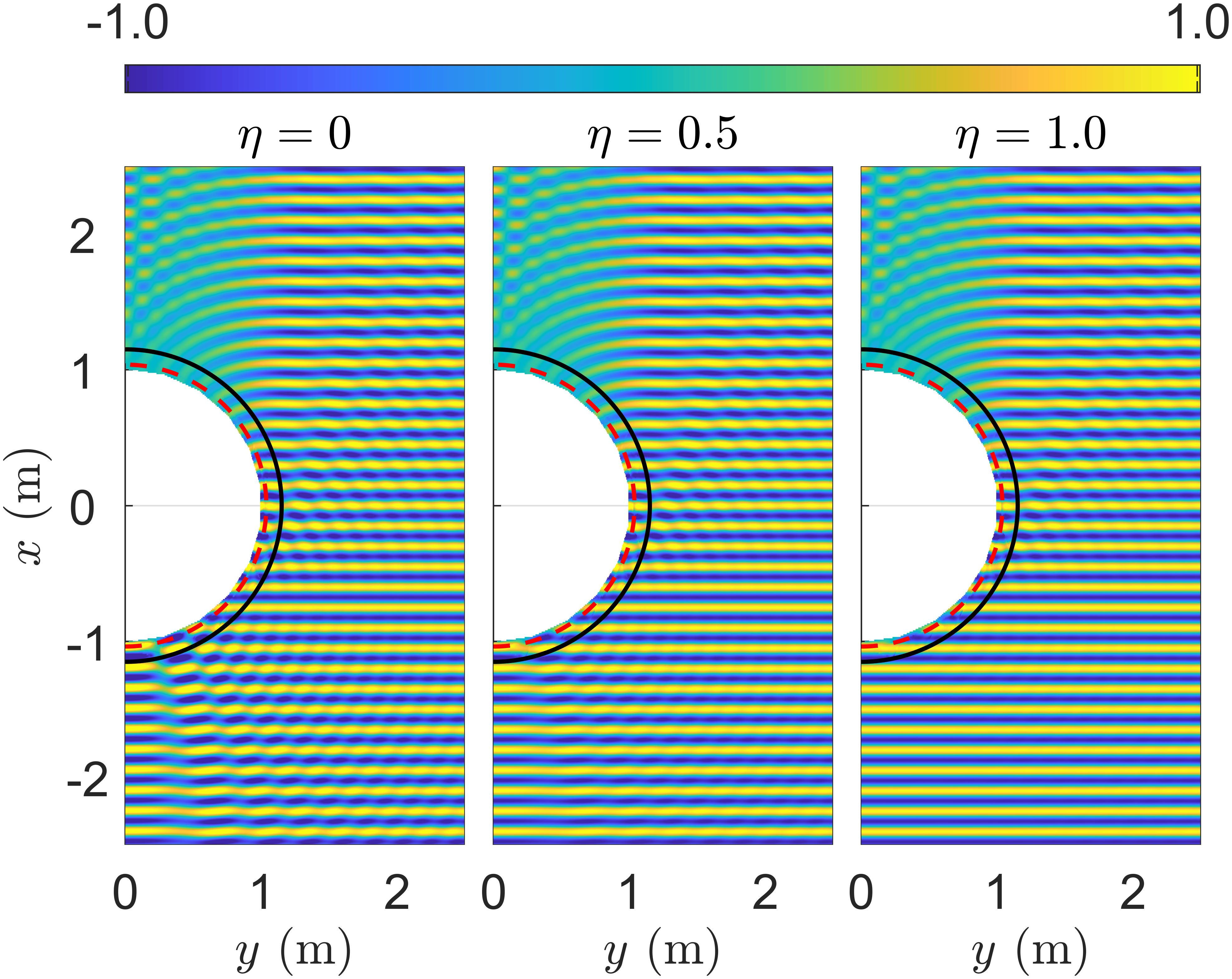}
    \caption{Snapshots of the total electric field $\Re(E_y)$ for the three choices of resistive sheet parameter $\eta$ (left: $\eta=0$, middle: $\eta=0.5$ and right: $\eta=1$). The dashed red line marks the resistive sheet interface and the solid black line the artificial interface used to introduce the incident wave. Backscattering appears less for $\eta=1$ compared to $\eta=0.5$ or the pure PEC sphere.}
    \label{fig:sheet_sphere_snap}
\end{figure}

\begin{figure}[!htb]
    \centering
    \includegraphics[width=0.6\textwidth]{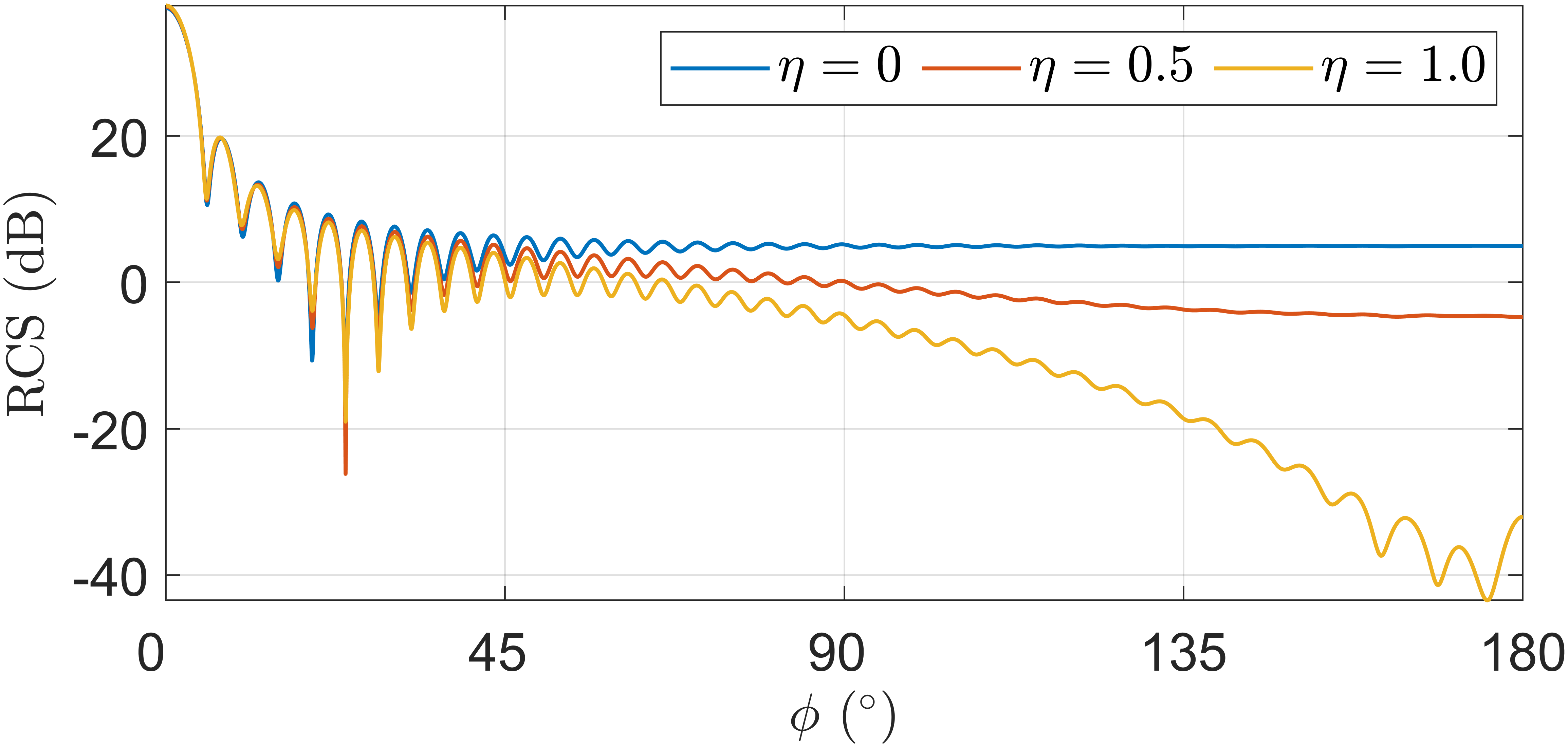} 
\caption{Bistatic RCS for the PEC sphere and resistive sheet example at 2 GHz.  The decreased backscattering due to the resitive sheet covering the sphere at $\eta=1$ is clearly seen.}
    \label{fig:sheet_sphere_mie}
\end{figure}

In Fig. \ref{fig:sheet_sphere_snap}, the total field component $\Re(E_y)$ on the plane $z=0$ is shown for three choices of $\eta$. We can no longer expect invisibilty since the screen is curved, but it is evident that backscattering is decreasing as  $\eta$ increases. This is seen more clearly in Fig.~\ref{fig:sheet_sphere_mie} where we show the RCS in each case.

\subsection{Heterogeneous models}\label{sec:het}

\subsubsection{A dielectric sphere}\label{penet_sphere}

In this experiment, a penetrable sphere with a radius of 1 meter is centered at the origin inside the cube $[-1-15\lambda_0, 1+15\lambda_0]^3$, where $\lambda_0$ is the wavelength in vacuum. For the penetrable sphere, we assume $\epsr = 1.5+0.5i$ and $\mur = 1$, while we select vacuum parameters in other domains. The frequency of the incident field is $f = 2$ GHz and a PML with thickness of $5\lambda_0$ is used on each side of the cube. An artificial spherical boundary with radius $1+\lambda_0$ is used to separate a scattered field region outside and
a total field region inside this surface.  The scattered-total field formulation in Section~\ref{stf} is used to introduce a source on the artificial boundary.

We utilised curved elements and a mesh size parameter of $\hs=3\lambda_s$, where $\lambda_s$ is a measure of the wavelength in the penetrable sphere computed using $\kappa_{\mathrm{abs}}$ defined in Appendix \ref{sec:app}. Figure \ref{fig:penet_sphere_grid} shows cross-sections of the computational grids.  To solve this problem, we used the supercomputer Puhti with a total of 7 computing nodes and 40 CPU units per node. More detailed information on the computations, including information on the grids and CPU time, is given in Table \ref{tab:penet}.

\begin{figure}[!htb]
    \centering
    \includegraphics[width=0.6\textwidth]{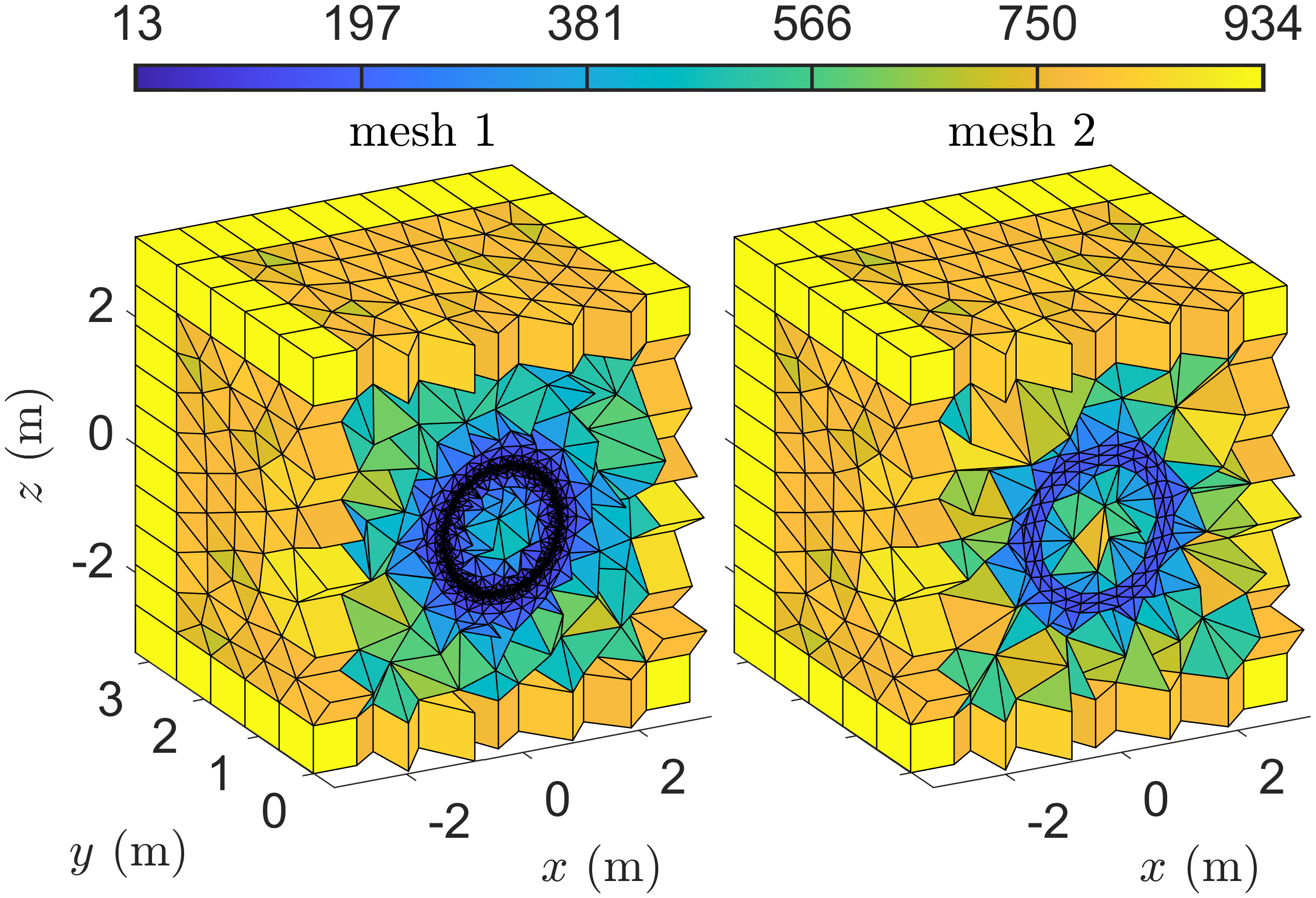} 
\caption{Cross-section of the computational grid for the penetrable sphere case. }
    \label{fig:penet_sphere_grid}
\end{figure}

\begin{figure}[!htb]
    \centering
    \includegraphics[width=0.55\textwidth]{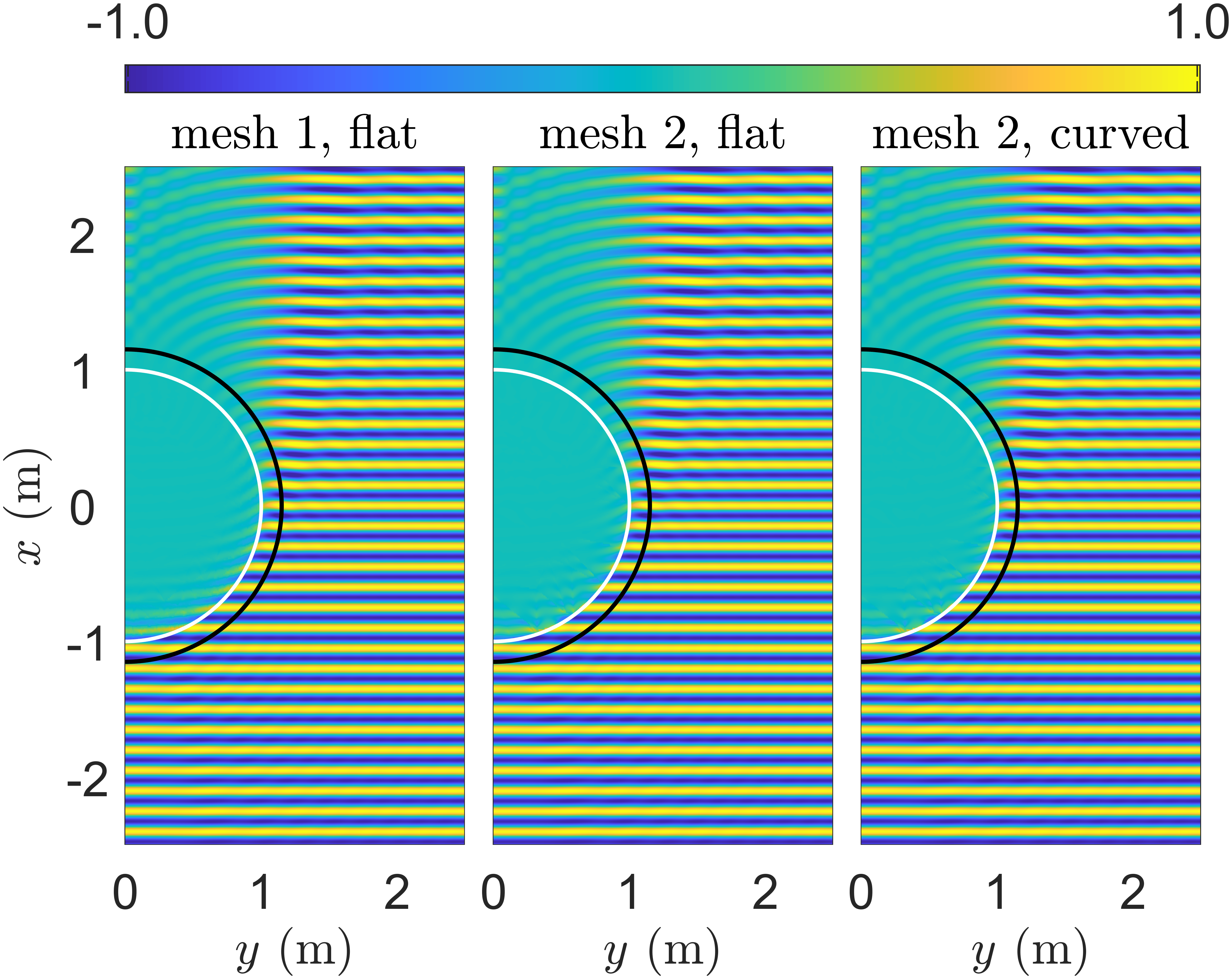}
    \caption{Snapshots of the electric field $\Re(E_y)$ on the $z=0$ plane for the penetrable sphere. The solid white line shows the material interface and the solid black line marks the artificial interface used to introduce the incident wave.}
    \label{fig:penet_sphere_snap}
\end{figure}

\begin{figure}[!htb]
    \centering
    \includegraphics[width=0.6\textwidth]{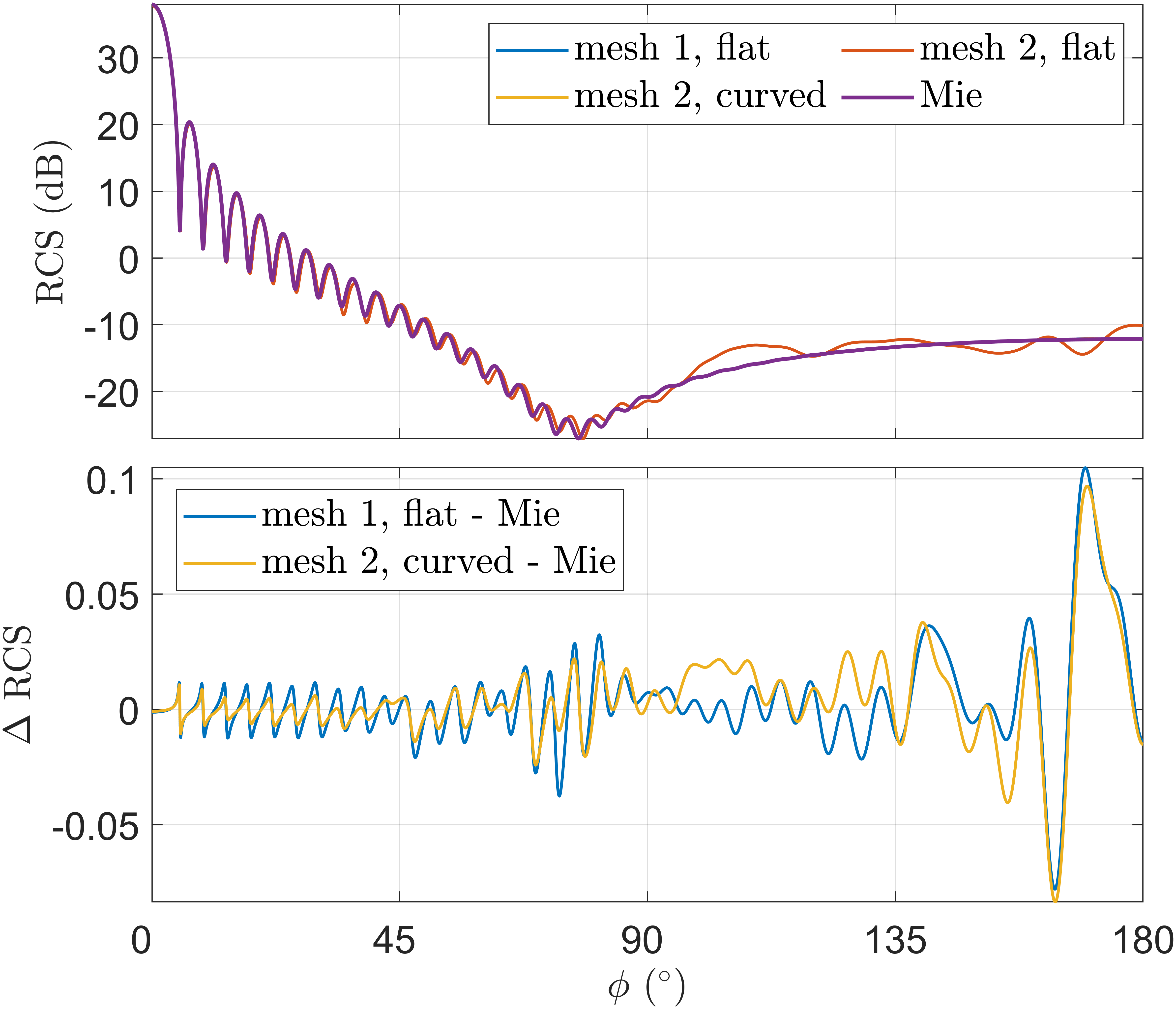} 
\caption{Bistatic RCS for the penetrable sphere at 2 GHz, otherwise the layout is the same caption as for Fig. \ref{fig:pec_sphere_mie}.}
    \label{fig:penet_sphere_mie}
\end{figure}

In Fig. \ref{fig:penet_sphere_snap}, we show the total field component $\Re(E_y)$ on the $z=0$-plane for the two meshes. As expected (see also Table~\ref{tab:penet}), \emph{mesh 2} with curved face
elements produces a solution comparable to \emph{mesh 1}, but faster.  The accuracy of the far field pattern is compared to the Mie series solution in Fig. \ref{fig:penet_sphere_mie}, and again demonstrates that \emph{mesh 1} and \emph{mesh 2} with curved faces give comparable results.

\begin{table}[!htb]
    \caption{Details of the computational results for the penetrable sphere. See Table~\ref{not} for definitions of the reported quantities.}
    \centering
    \resizebox{\textwidth}{!}{
      \begin{tabular}{cc|ccccccc|ccc}
        mesh id & surface & $\Nt$ & $\Nw$ & $\Nh$ & $\Nv$ & $h_{\min}$ (cm) & $h_{\max}$ (m) & $\Ndof$ & $N_{{\mathrm{iter}}}^{{\mathrm{BiCG}}}$ & L2-error (\%) & CPU-time (s)\\
        \hline
        1 & flat & 337,550 & 936 & 104 & 59,053 & 0.65 & 1.30 & 21,633,354 & 187 & 0.14 & 346\\
        2 & flat & 7,100 & 936 & 104 & 2,161 & 7.80 & 1.42 & 4,640,348 & 90 & 7.85 & 162\\
        2 & curved & 7,100 & 936 & 104 & 2,161 & 7.80 & 1.42 & 4,640,348 & 88 & 0.14 & 234
      \end{tabular}
      }
    \label{tab:penet}
\end{table}

For a comparison of \emph{ParMax} with edge finite elements in the case of the dielectric sphere, see Appendix~\ref{sec:app_comp}.

\subsubsection{A plasma sphere}

In this second experiment for penetrable objects, we assume a penetrable sphere with a radius of 0.25 meter is centred at the origin of a cube $[-0.25-15\lambda_0, 0.25+15\lambda_0]^3$.  The material in the sphere is assumed to be a plasma modelled by setting $\epsr = -1.5+0.5i$ and $\mur = 1$. The frequency of the incident field is $f = 2$ GHz and a perfectly matched layer with thickness of $5\lambda_0$ is used inside each side of the cube. The source is introduced on an artificial spherical surface  of  radius $0.25+\lambda_0$.

We utilise curved elements and a mesh size parameter  $\hs=3\lambda_s$, where $\lambda_s$ denotes the wavelength in the penetrable sphere, on the material discontinuity surface. The grid comprises of 2,959 tetrahedral elements that cover the main domain of interest. Furthermore, the surrounding PML layer is discretized using a combination of 92 hexahedral and 768 wedge elements. The entire grid has 1,304 vertices, with $h_{\min} = 5.88$ cm and $h_{\max} = 1.30$ m. Here, the degrees of freedom number $\Ndof$ is 2,234,416.

\begin{figure}[!htb]
    \centering
    \includegraphics[width=0.55\textwidth]{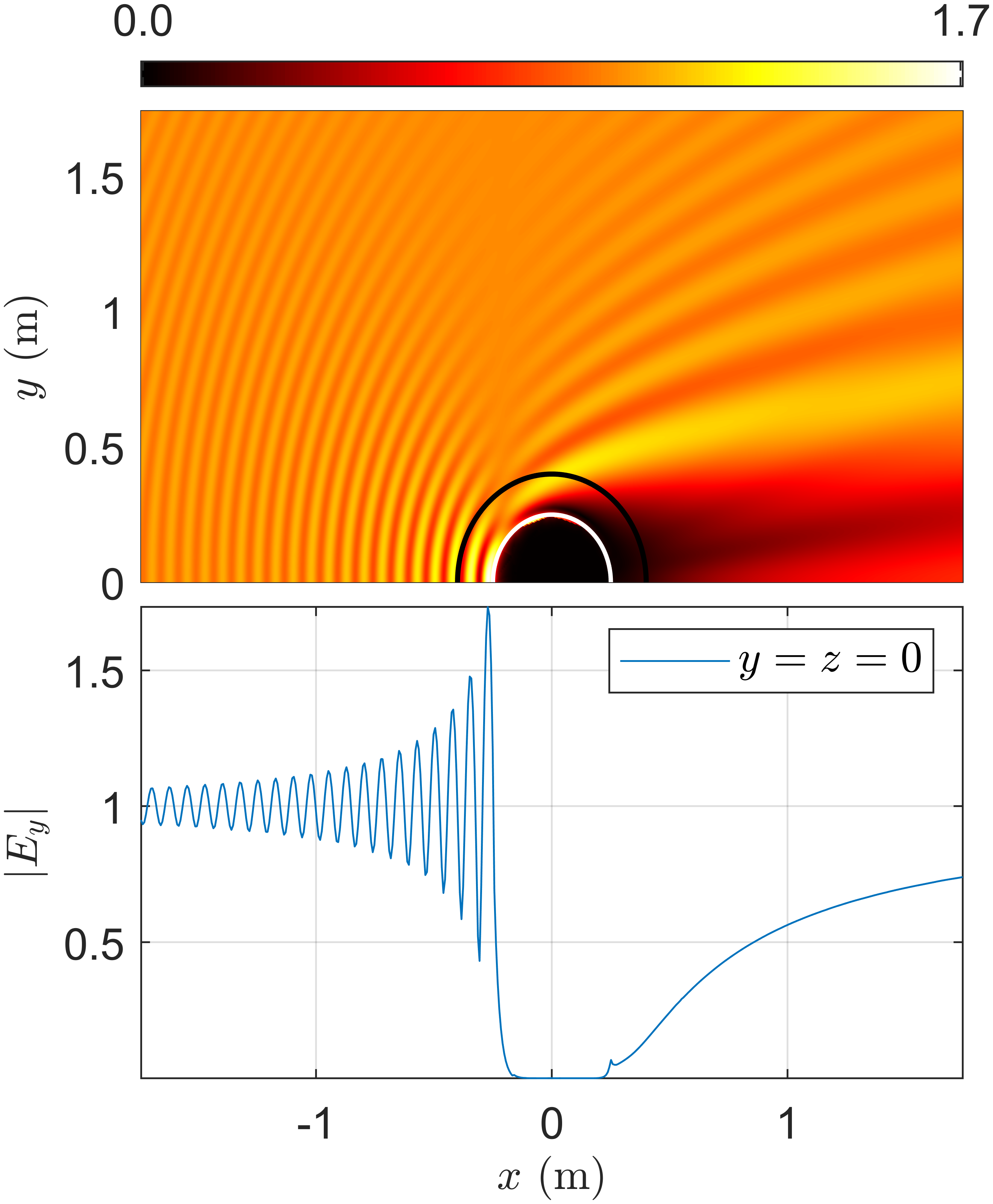}
    \caption{Snapshots of the total electric field. Top panel: a contour plot of $|E_y|$ in the plane $z=0$. The solid white line shows the interface between vacuum and the plasma sphere. The solid black line marks the interface used to introduce the incident wave. Bottom panel: a plot of $|E_y|$ along the line $y=z=0$.}
    \label{fig:plasma_sphere_snap}
\end{figure}

\begin{figure}[!htb]
    \centering
    \includegraphics[width=0.6\textwidth]{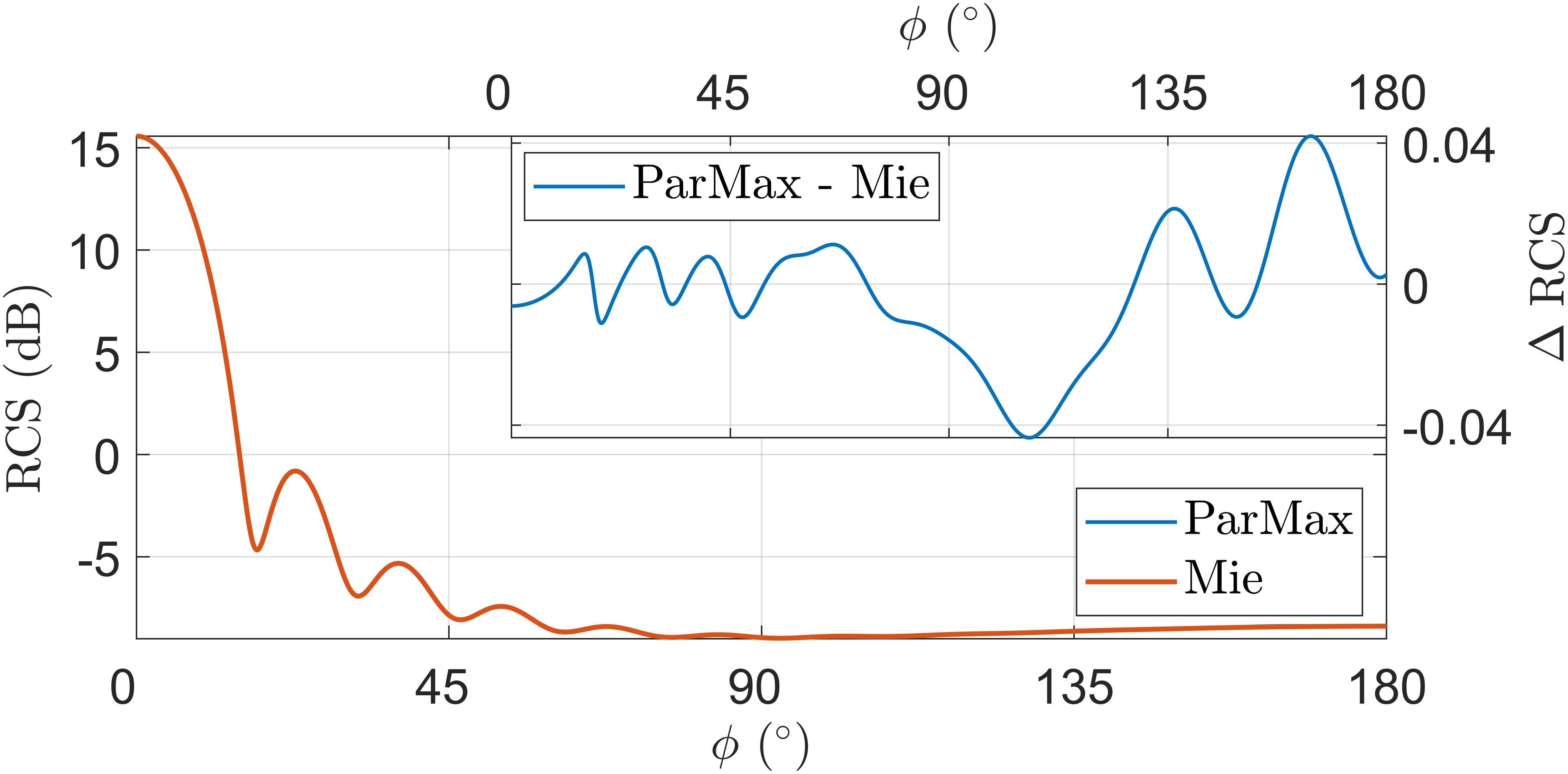} 
\caption{Bistatic RCS for the plasma sphere at 2.0 GHz, comparing the computed RCS and Mie result. Top right corner shows the difference between the numerical and Mie series.}
    \label{fig:plasma_sphere_mie}
\end{figure}

We used the supercomputer Puhti with a total of 5 computing nodes and 10 CPU units per node to solve the problem. It took 427 seconds from each CPU unit to build the system matrices $C$ and $D$ and then reach the solution with 103 bi-conjugate iterations. Snapshots of the total field are shown in Fig.~\ref{fig:plasma_sphere_snap}.  In Fig.~\ref{fig:plasma_sphere_mie} we show the computed RCS which shows remarkable agreement with the Mie series solution.

\section{Conclusion}\label{sec:conc}

In this study, we explored and extended the Ultra-Weak Variational Formulation (UWVF) applied to the time-harmonic Maxwell's equations. Our research findings led to important contributions that enhance the efficiency and applicability of the UWVF method for solving electromagnetic wave problems.

The paper shows a series of numerical examples validating the effectiveness of the new enhancements. Scattering problems from PEC objects were considered, highlighting the benefits of curved elements, different element types, and the low-memory version of the software. The applicability was further demonstrated through simulations of scattering from an full size aircraft, emphasising its potential for real-world industrial scenarios.

This paper has not only extended the capabilities of the UWVF for electromagnetic wave problems but has also provided a comprehensive set of numerical results to underline the practical significance of these advancements. The integration of curved elements, various element shapes, and resistive sheets collectively contribute to the method's robustness and utility, making it a valuable tool for addressing complex electromagnetic problems.

The choice of mesh size for the surface triangulation needs further study in the case of large imaginary part of $\epsr$ or in regions of high curvature but these issues are beyond the scope of the paper. Related to this, an important direction for further work would be to refine the 
heuristics for choosing the number and direction of plane waves on
each element. This is particularly needed for elements that might have
large aspect ratio such as elongated wedges.

\section*{Acknowledgments}

The research of P. M. is partially supported by the US AFOSR under grant number FA9550-23-1-0256.  The research of T. L. is partially supported by the Academy of Finland (the Finnish Center of Excellence of Inverse Modeling and Imaging) and the Academy of Finland project 321761. The authors also wish to acknowledge the CSC – IT Center for Science, Finland, for generously sharing their computational resources.

\appendices
\section{Choice of basis}\label{sec:app}

\begin{figure}[!htb]
    \centering
    \includegraphics[height=0.215\textwidth]{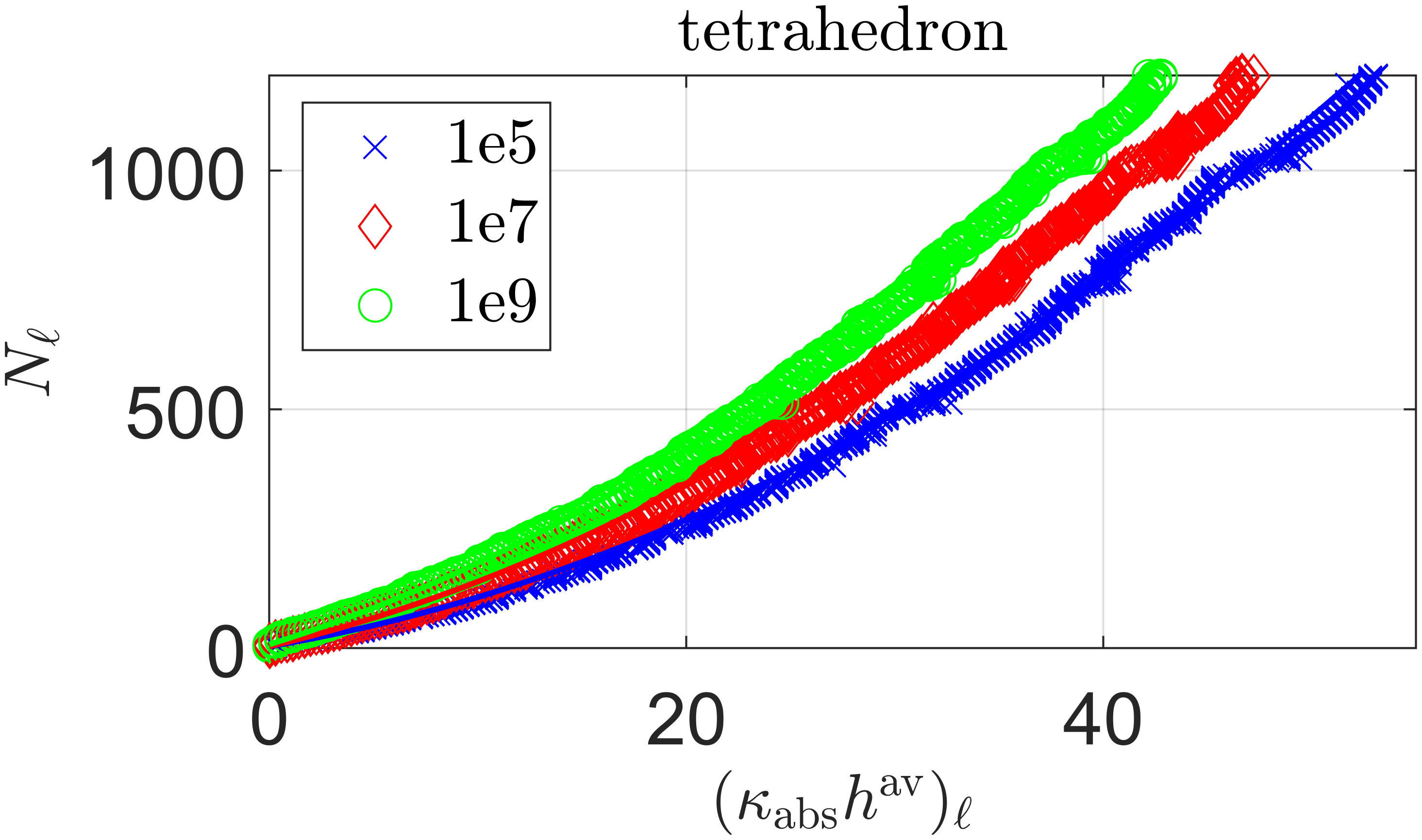} 
    \includegraphics[height=0.215\textwidth]{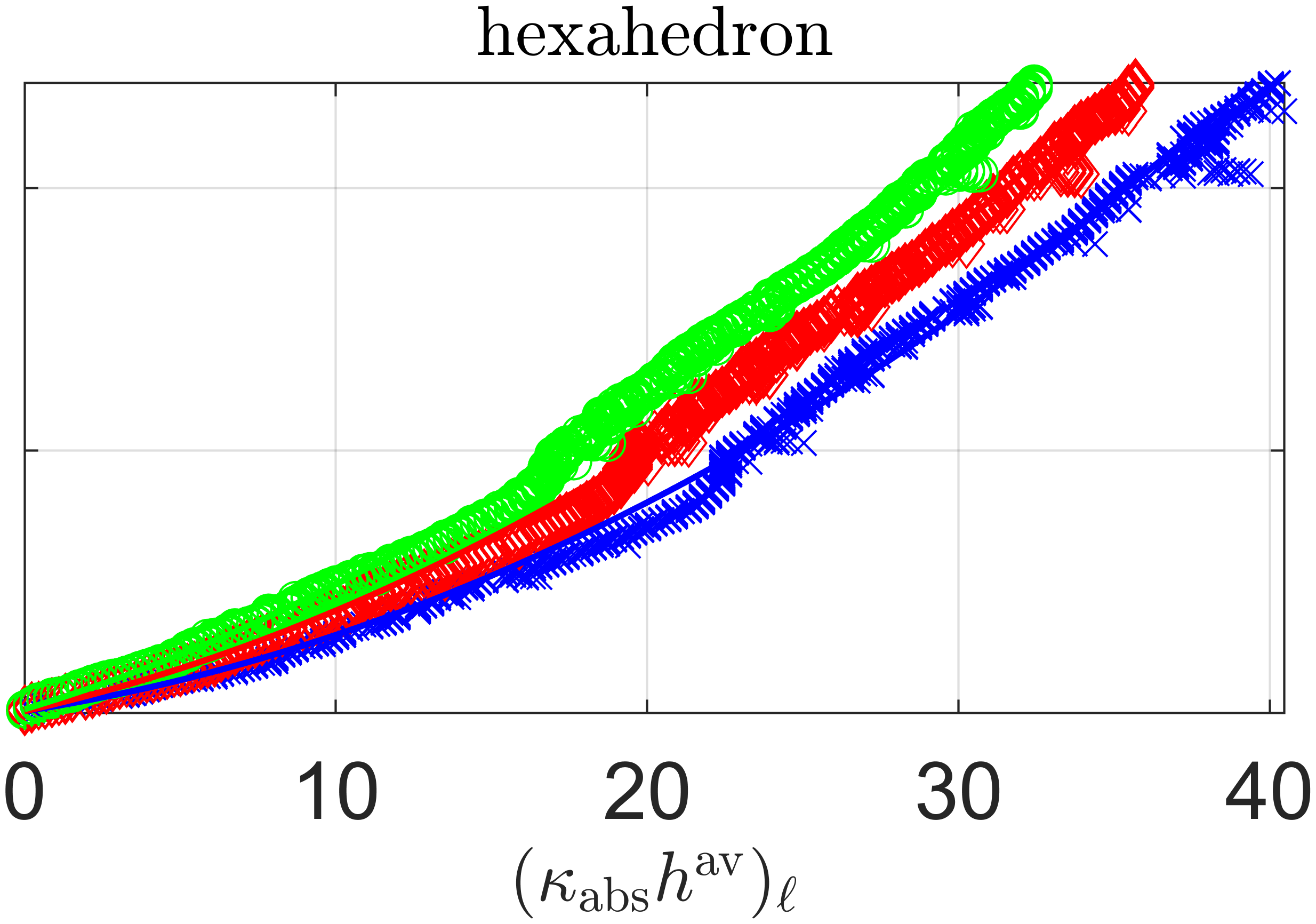} 
    \includegraphics[height=0.215\textwidth]{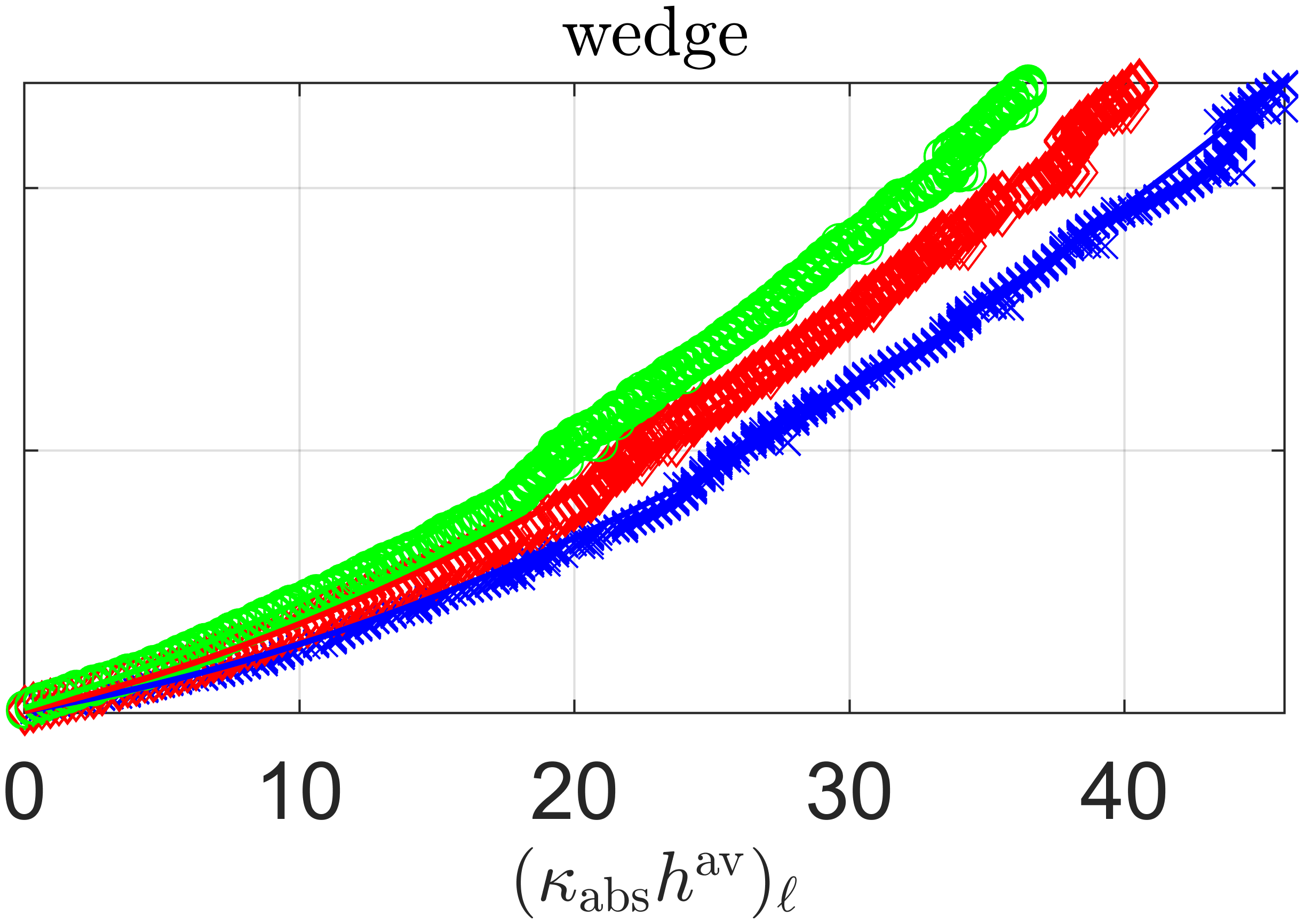} 
\caption{The number of basis functions $N_{\ell}$ as a function of $(\kappa_{\mathrm{abs}}h^{\mathrm{av}})_{\ell}$ when the basis dimension is chosen by constraining the maximum
condition number of $D_{\ell}$.}
    \label{fig:pw_selection}
\end{figure}

Because we have used new element types and larger numbers of directions in this paper compared to \cite{Huttunen2007}, we need new heuristics for choosing the number of plane wave directions on a particular element.  We use the same technique as in \cite{Huttunen2007}. Computing on the reference element, in Fig. \ref{fig:pw_selection}, the number of plane waves $N_{\ell}$ is plotted as a function of $(\kappa_{\mathrm{abs}}h^{\mathrm{av}})_{\ell}$, when the maximum condition number of the matrix blocks of $D_{\ell}$ is limited by the tolerances $10^5$, $10^7$, and $10^9$.
Here
\begin{displaymath}
    \kappa_{\mathrm{abs}} = \omega\left|\sqrt{\epsr \mur}\right| 
\end{displaymath}
and the  element size parameter $h^{\mathrm{av}}$ is defined as the mean distance of the element's vertices from its centroid. 

This data is fitted by a quadratic polynomial function with the constraint that the polynomial gives at least 4 directions even on the finest grid:
\begin{equation}
\label{eq:Nell}
N_{\ell} = \left \lceil{a\left(\kappa_{\mathrm{abs}}h^{\mathrm{av}}\right)_{\ell}^2 + b\left(\kappa_{\mathrm{abs}}h^{\mathrm{av}}\right)_{\ell} + c}\right \rceil.
\end{equation}

Results of this fitting are shown in Table~\ref{tab:values}. Although only computed for one element shape, these polynomials are used to set the number of directions for any given mesh.  Generally a higher tolerance on the condition number results in more directions per element so greater accuracy, but too high a condition number slows BiCGstab unacceptably.

\begin{table}[!htb]
    \caption{Parameters for the basis polynomials (see (\ref{eq:Nell})).}
    \centering
    \resizebox{\textwidth}{!}{
      \begin{tabular}{c|ccc|ccc|ccc}
        & & tetrahedron & & & hexahedron & & & wedge  & \\
        max(cond($D_k$)) & $a$ & $b$ & $c$ & $a$ & $b$ & $c$ & $a$ & $b$ & $c$\\
        \hline
        1e5 & 0.2972 & 7.3336 & 4.0000 & 0.5365 & 9.1369 & 4.0000 & 0.3752 & 9.0041 & 4.0000\\
        1e7 & 0.3305 & 10.2707 & 4.0000 & 0.5803 & 13.3338 & 4.0000 & 0.4325 & 12.2717 & 4.0000\\
        1e9 & 0.3430 & 13.6221 & 8.1296 & 0.5967 & 17.7977 & 7.7490 & 0.4704 & 15.6097 & 9.9414\\
      \end{tabular}
      }
    \label{tab:values}
\end{table}

\section{Comparison to edge elements}\label{sec:app_comp}

There are many ways to solve the scattering example problems in this
paper. To provide a comparison to a more familiar method, we now use
\emph{ParMax} and the edge finite element method to compute scattering
from a dielectric sphere.  The exact Mie series is used to assess
accuracy at several frequencies. 

We have used the open source Netgen/NGSolve finite element
package~\cite{netgen}.  This is a library of highly optimised C++
function supporting multithreading that can be used via a Python
interface (called NGSpy).  It includes the Netgen mesh generator.  We
use $p$th order N\'ed\'elec edge elements of full polynomial degree on
a tetrahedral mesh with a rectilinear PML and an outer impedance
boundary condition.  A direct sparse matrix solver supplied with NGSpy
was used to solve the linear system.  This limited the overall size of
the problem we could solve to approximately 1.2 million degrees of
freedom.

In order to make the application of the FEM possible, the
computational domain described in Section~\ref{sec:het} is not used
here. Instead, the computational domain is the cube
$[-1-(1+d_{\mathrm{pml}})\lambda_0,1+(1+d_{\mathrm{pml}})\lambda_0]^3$. A
PML of relative width $d_{\mathrm{pml}}=0.3$ with FEM and
$d_{\mathrm{pml}}=1$ with \emph{ParMax} is added within this cube. For
the FEM, we requested NGSpy to use the mesh size
$\lambda_0(2p+1)/(4\pi)$ outside the scatterer and increased $p$
gradually from $p=3$ to $p=8$ as $\lambda_0$ decreased to maintain an
accuracy of roughly 1\% in the computed RCS. Inside the scatterer the
mesh size is reduced by a factor of 0.4.

\begin{figure}[!hbt]
    \centering
    \includegraphics[width=0.6\textwidth]{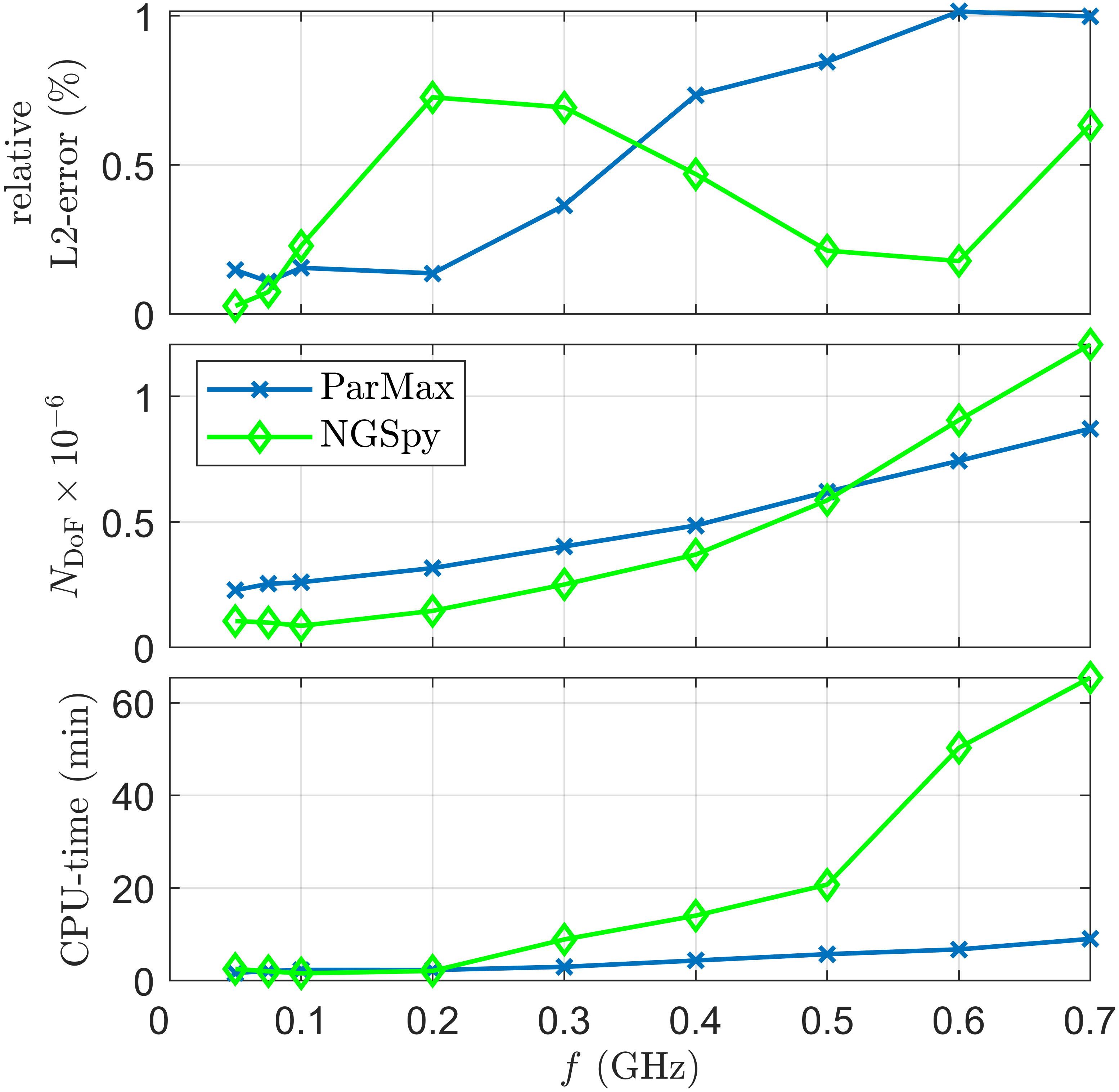}
\caption{Comparison results between NGSpy and \emph{ParMax}.}
    \label{fig:comparison}
\end{figure}

The computations were run on an Intel Xeon Gold 6138 CPU \@ 2.00GHz
with 20 cores, each having two threads per core, and 394GB of
memory. Since NGSpy ran multithreaded computations, we ran
\emph{ParMax} with MPI and 20 cores. Elapsed time is reported for
NGSpy and \emph{ParMax}, including assembly, matrix factorization,
solving, and computation of the far-field pattern. The error reported
is the relative percentage error in the RCS compared to the Mie
solution. The results are shown in Fig.~\ref{fig:comparison}.

Figure~\ref{fig:comparison} (top) shows that the relative L2-error
remains comparable between NGSpy and \emph{ParMax} across the entire
frequency spectrum used here. As illustrated in
Fig.~\ref{fig:comparison} (middle), the number of degrees of freedom
required by NGSpy increases more rapidly than that for
\emph{ParMax}. Likewise, the comparison in Fig.~\ref{fig:comparison}
(bottom) reveals that the CPU-time required by NGSpy increases more
significantly compared to \emph{ParMax}. This could perhaps be
ameliorated using an auxiliary space preconditioned iterative
scheme~\cite{HiptmairXu} rather than a direct factorization, but this
is not available to us in NGSpy.


\begin{thebibliography}{10}

\bibitem{CSC}
{CSC – IT C}enter for {S}cience {L}td, {C}omputing environments.
\newblock \url{https://docs.csc.fi/computing/available-systems/}.
\newblock Accessed: Aug. 1, 2023.

\bibitem{badics14}
Z.~Badics.
\newblock Trefftz-{D}iscontinuous {G}alerkin and finite element multi-solver
  technique for modeling time-harmonic {EM} problems with high-conductivity
  regions.
\newblock {\em IEEE Trans. on Magnetics.}, 50(2):401--404, 2014.

\bibitem{Berenger}
J.-P. B\'erenger.
\newblock A {P}erfectly {M}atched {L}ayer for the absorption of electromagnetic
  waves.
\newblock {\em J. Comput. Phys.}, 114(2):185--200, 1994.

\bibitem{buf07}
A.~Buffa and P.~Monk.
\newblock Error estimates for the {U}ltra {W}eak {V}ariational {F}ormulation of
  the {H}elmholtz equation.
\newblock {\em ESAIM: Mathematical Modeling and Numerical Analysis},
  42:925--40, 2008.

\bibitem{cessenat_phd}
O.~Cessenat.
\newblock {\em Application d'une nouvelle formulation variationnelle aux
  \'{e}quations d'ondes harmoniques. {P}robl\`{e}mes de {H}elmholtz 2{D} et de
  {M}axwell 3{D}.}
\newblock PhD thesis, Universit\'{e} Paris IX Dauphine, 1996.

\bibitem{cessenat03}
O.~Cessenat and B.~Despr\'{e}s.
\newblock Using plane waves as base functions for solving time harmonic
  equations with the {U}ltra {W}eak {V}ariational {F}ormulation.
\newblock {\em J. Comput. Acoust.}, 11:227--38, 2003.

\bibitem{git07}
C.~Gittelson, R.~Hiptmair, and I.~Perugia.
\newblock Plane wave discontinuous {G}alerkin methods.
\newblock {\em ESAIM: Mathematical Modeling and Numerical Analysis},
  43:297--331, 2009.

\bibitem{hafner85}
C.~Hafner.
\newblock {MMP} calculations of guided waves.
\newblock {\em IEEE Trans. On Magn.}, 21, 1985.

\bibitem{hardin_16}
D.~Hardin, T.~Michaels, and E.~Saff.
\newblock A comparison of popular point configurations on ${S}^2$.
\newblock {\em Dolomites Research Notes on Approximation}, 9:16--19, 2016.

\bibitem{moi11}
R.~Hiptmair, A.~Moiola, and I.~Perugia.
\newblock Stability results for the time-harmonic {M}axwell equations with
  impedance boundary conditions.
\newblock {\em Math. Meth. Appl. Sci.}, 21:2263--87, 2011.

\bibitem{HMP13}
R.~Hiptmair, A.~Moiola, and I.~Perugia.
\newblock Error analysis of {T}refftz-discontinuous {G}alerkin methods for the
  time-harmonic {M}axwell equations.
\newblock {\em Math. Comput.}, 82:247--268, 2013.

\bibitem{HiptmairXu}
R.~Hiptmair and J.~Xu.
\newblock Nodal auxiliary space preconditioning in ${H}({\rm curl})$ and
  ${H}({\rm div})$ spaces.
\newblock {\em SIAM J. Numer. Anal.}, 45:2483--2509, 2007.

\bibitem{Huttunen2007}
T.~Huttunen, M.~Malinen, and P.~Monk.
\newblock Solving {M}axwell's equations using the ultra weak variational
  formulation.
\newblock {\em J. Comput. Phys.}, 223:731--758, 2007.

\bibitem{IM-sylvand-21}
L.~Imbert-G\'erard and G.~Sylvand.
\newblock A roadmap for {G}eneralized {P}lane {W}aves and their interpolation
  properties.
\newblock {\em Numer. Math.}, 149:87–137, 2021.

\bibitem{jin_volakis}
J.~Jin, J.~Volakis, C.~Yu, and A.~Woo.
\newblock Modeling of resistive sheets in finite element solutions.
\newblock {\em IEEE Trans. Antennas Propag.}, 40:727--731, 1992.

\bibitem{Karypis_metis98}
G.~Karypis and V.~Kumar.
\newblock A fast and high quality multilevel scheme for partitioning irregular
  graphs.
\newblock {\em SIAM J. Sci. Comput.}, 20(1):359--392, 1998.

\bibitem{PS_22}
C.~Lehrenfeld and P.~Stocker.
\newblock Embedded {T}refftz discontinuous {G}alerkin methods.
\newblock {\em Int. J. Numer. Methods Eng.}, 124(17):3637--3661, 2023.

\bibitem{macphie03}
R.~H. MacPhie and K.-L. Wu.
\newblock A plane wave expansion of spherical wave functions for modal analysis
  of guided wave structures and scatterers.
\newblock {\em IEEE Trans. Antennas Propag.}, 51:2801--2805, 2003.

\bibitem{Monk03}
P.~Monk.
\newblock {\em Finite Element Methods for {M}axwell's Equations}.
\newblock Oxford University Press, Oxford, 2003.

\bibitem{sirdey22}
S.~Pernet, M.~Sirdey, and S.~Tordeux.
\newblock Ultra-weak variational formulation for heterogeneous {M}axwell
  problem in the context of high performance computing.
\newblock Available as \url{https://hal.science/hal-03642116} from the Hal
  preprint server., 2022.

\bibitem{netgen}
J.~Sch\"oberl.
\newblock Netgen/{NGS}olve.
\newblock Downloaded from \url{https://ngsolve.org}, 2023.

\bibitem{Senior}
T.~Senior.
\newblock Combined resistive and conductive sheets.
\newblock {\em IEEE Trans. Antennas Propag.}, AP-33:577--579, 1985.

\bibitem{NGSTrefftz}
P.~Stocker.
\newblock {NGST}refftz: {A}dd-on to {NGS}olve for {T}refftz methods.
\newblock {\em J. Open Source Soft.}, 7, 2022.
\newblock Art. No. 4135.

\bibitem{Trefftz26}
E.~Trefftz.
\newblock Ein gegenst\"{u}ck zum {R}itz'schen verfahren.
\newblock In {\em Proc. 2nd Int. Congr. Appl. Mech.}, pages 131--137, Zurich,
  1926.

\end{thebibliography}

\end{document}